\documentclass[12pt,reqno]{amsart}
\usepackage{amsmath , amssymb}
\usepackage[mathscr]{eucal}

\newdimen\unit\newdimen\psep\newcount\nd\newcount\ndx\newbox\dotb\newbox\ptbox
\newdimen\dx\newdimen\dy\newdimen\dxx\newdimen\dyy\newdimen\hgt
\newdimen\xoff\newdimen\yoff
\newcommand\clap[1]{\hbox to 0pt{\hss{#1}\hss}}
\newcommand\vdisk[1]{{\font\dotf=cmr10 scaled #1\dotf.}}
\newcommand\varline[2]{\setbox\dotb\hbox{\vdisk{#1}}\xoff=-.5\wd\dotb
\wd\dotb=0pt\yoff=-.5\ht\dotb\psep=#2\ht\dotb}
\newcommand\varpt[1]{\setbox\ptbox\clap{\vdisk{#1}}\setbox\ptbox
\hbox{\raise-.5\ht\ptbox\box\ptbox}}
\newcommand\cpt{\copy\ptbox}
\newcommand\point[3]{\rlap{\kern#1\unit\raise#2\unit\hbox{#3}}}
\newcommand\setnd[4]{\dx=#3\unit\advance\dx-#1\unit\divide\dx by\psep
\dy=#4\unit\advance\dy-#2\unit\divide\dy by\psep \multiply\dx
by\dx\multiply\dy by\dy\advance\dx\dy\nd=1\advance\dx-1sp
\loop\ifnum\dx>0\advance\dx-\nd sp\advance\nd1\advance\dx-\nd
sp\repeat}

\newcommand\dline[5]{{\nd=#5\hgt=#2\unit\dx=#3\unit\advance\dx-#1\unit
\divide\dx by\nd\dy=#4\unit\advance\dy-#2\unit\divide\dy by\nd
\advance\hgt\yoff\rlap{\kern#1\unit\kern\xoff\loop\ifnum\nd>1\advance\nd-1
\advance\hgt\dy\kern\dx\raise\hgt\copy\dotb\repeat}}}

\newcommand\qellip[4]{{\setnd{0}{0}{#3}{#4}\dx=\unit\dy=0pt\raise\yoff\rlap{%
\kern#1\unit\kern\xoff\raise#2\unit\hbox{\loop\ifnum\dx>0\rlap{\kern#3\dx
\raise#4\dy\copy\dotb}\hgt=\dx\divide\hgt
by\nd\advance\dy\hgt\hgt=\dy \divide\hgt
by\nd\advance\dx-\hgt\repeat\rlap{\raise#4\dy\copy\dotb}}}}}
\newcommand\bez[6]{{\setnd{#1}{#2}{#3}{#4}\ndx=\nd\setnd{#3}{#4}{#5}{#6}
\ifnum\ndx>\nd\nd=\ndx\fi\dx=#3\unit\advance\dx-#1\unit\dy=#4\unit
\advance\dy-#2\unit\dxx=#5\unit\advance\dxx-#1\unit\dyy=#6\unit\advance
\dyy-#2\unit\advance\dxx-2\dx\advance\dyy-2\dy\divide\dxx
by\nd\divide\dyy
by\nd\advance\dx.25\dxx\advance\dy.25\dyy\divide\dx
by\nd\divide\dy by\nd \multiply\nd
by2\dx=100\dx\dy=100\dy\dxx=100\dxx\dyy=100\dyy\divide\dxx by\nd
\divide\dyy
by\nd\hgt=#2\unit\raise\yoff\rlap{\kern#1\unit\kern\xoff
\raise\hgt\copy\dotb\loop\ifnum\nd>0\advance\nd-1\advance\hgt0.01\dy
\kern0.01\dx\raise\hgt\copy\dotb\advance\dx\dxx\advance\dy\dyy\repeat}}}
\newcommand\ptu[3]{\point{#1}{#2}{\cpt\raise1ex\clap{$\scriptstyle{#3}$}}}
\newcommand\ptd[3]{\point{#1}{#2}{\cpt\raise-1.8ex\clap{$\scriptstyle{#3}$}}}
\newcommand\ptr[3]{\point{#1}{#2}{\cpt\raise-.4ex\rlap{$\
\scriptstyle{#3}$}}}
\newcommand\ptl[3]{\point{#1}{#2}{\cpt\raise-.4ex\llap{$\scriptstyle{#3}\
$}}}
\newcommand\ptlu[3]{\point{#1}{#2}{\raise.8ex\clap{$\scriptstyle{#3}$}}}
\newcommand\ptld[3]{\point{#1}{#2}{\raise-1.6ex\clap{$\scriptstyle{#3}$}}}
\newcommand\ptlr[3]{\point{#1}{#2}{\raise-.4ex\rlap{$\,\scriptstyle{#3}$}}}
\newcommand\ptll[3]{\point{#1}{#2}{\raise-.4ex\llap{$\scriptstyle{#3}\,$}}}

\newcommand\thnline{\varline{400}{.6}}

\varpt{2500}\thnline\unit=2em

\newtheorem{thm}{Theorem}
\newtheorem{conj}{Conjecture}

\newtheorem*{mader}{Mader's Theorem}

\newtheorem{qu}{Question}
\newtheorem{prob}{Problem}
\newtheorem{lemma}[thm]{Lemma}
\newtheorem{cor}[thm]{Corollary}

\newtheorem{obs}{Observation}
\theoremstyle{definition}\newtheorem{rmk}{Remark}
\theoremstyle{definition}
\newcommand{\ds}{\displaystyle}

\newcommand{\ol}{\overline}
\def\N{\mathbb{N}}
\def\RR{\mathbb{R}}
\def\P{\mathcal{P}}
\def\S{\mathcal{S}}
\def\T{\mathcal{T}}
\def\eps{\varepsilon}
\def\le{\leqslant}
\def\ge{\geqslant}

\begin{document}

\title[Highly connected multicoloured subgraphs]{Highly connected multicoloured subgraphs of multicoloured graphs}

\author{Henry Liu}\address{Department of Mathematical Sciences\\ The University of Memphis\\ Memphis, TN, 38152} \email{henryliu@memphis.edu}

\author{Robert Morris}\address{Department of Mathematical Sciences\\ The University of Memphis\\ Memphis, TN, 38152} \email{rdmorrs1@memphis.edu}

\author{Noah Prince}\address{Department of Mathematics\\ University of Illinois\\ 1409 W. Green Street\\ Urbana, IL 61801} \email{nprince@math.uiuc.edu}
\thanks{The second and third authors were partially supported during this research by Van Vleet Memorial Doctoral Fellowships.}

\maketitle

\begin{abstract}
Suppose the edges of the complete graph on $n$ vertices, $E(K_n)$, are coloured using $r$ colours; how large a $k$-connected subgraph are we guaranteed to find, which uses only at most $s$ of the colours? This question is due to Bollob\'as, and the case $s = 1$ was considered in~\cite{HNR1}. Here we shall consider the case $s \ge 2$, proving in particular that when $s = 2$ and $r+1$ is a power of 2 then the answer lies between $\frac{4n}{r+1} \: - \: 5kr(r + 2k + 1)$ and $\frac{4n}{r+1} \: + \: 4$, that if $r = 2s + 1$ then the answer lies between $(1 - 1/{r \choose s})n \: - \: 2{r \choose s}k$ and $(1 - 1/{r \choose s})n \: + \: 1$, and that phase transitions occur at $2s = r$ and $s = \Theta(\sqrt{r})$. We shall also mention some of the more glaring open problems relating to this question.
\end{abstract}

\section{Introduction}\label{HNR2intro}

In this paper we shall study the following extremal problem, due to Bollob\'as: given an $r$-colouring of the edges of the complete graph on $n$ vertices, $E(K_n)$, how large a $k$-connected subgraph can we find which uses at most $s$ of the colours? (Recall that a graph $G$ on $n \ge k+1$ vertices is said to be \textit{$k$-connected} if whenever at most $k-1$ vertices are removed from $G$, the remaining vertices are still connected by edges of $G$.) The case $s = 1$ of Bollob\'as' question was considered by the authors in a previous paper~\cite{HNR1}, and asymptotically tight bounds were obtained when $r-1$ is a prime power. In this paper we shall continue the investigations of \cite{HNR1} by considering the case $s \ge 2$, and in particular the case $s = 2$, and the cases $2s \approx r$ and $s = \Theta(\sqrt{r})$, where the function `jumps'. The majority of the problem remains wide open however, and so we shall also discuss some open problems and conjectures.

Let us begin by recalling the results and notation of \cite{HNR1}. We note that a gentler introduction into the problem is provided in that paper. Suppose we are given $n,r,s,k \in \N$, and a function $f: E(K_n) \to [r]$, i.e., an $r$-colouring of the edges of $K_n$. We assume always that $n \ge 2$. Given a subgraph $H$ of $K_n$, write $c_f(H)$ for the order of the image of $E(H)$ under $f$, i.e., $c_f(H) = |f(E(H))|$, the number of different colours with which $f$ colours $H$. Now, define $M(f,n,r,s,k) = \max\{|V(H)| : H \subset K_n$, $c_f(H) \le s\}$, the order of the largest $k$-connected subgraph of $K_n$ using at most $s$ colours from $[r]$. Finally, define $m(n,r,s,k) = \min_f\{M(f,n,r,s,k)\}$. Thus, the question of Bollob\'as asks for the determination of $m(n,r,s,k)$ for all values of the parameters. We shall state all of our results in terms of $m(n,r,s,k)$.

In \cite{HNR1} fairly tight bounds were given on the function $m(n,r,1,k)$. To be precise, it was shown that $m(n,2,1,k) = n - 2k + 2$ for every $n \ge 13k - 15$, that $$\frac{n-k+1}{2} \; \le \; m(n,3,1,k) \; \le \; \frac{n-k+1}{2} \: + \: 1$$ for every $n \ge 480k$, and more generally that $$\frac{n}{r-1} \: - \: 11(k^2 - k)r \; \le \; m(n,r,1,k) \; \le \; \frac{n-k+1}{r-1} \: + \: r$$ whenever $r - 1$ is a prime power.

In this paper we shall study the function $m(n,r,s,k)$ when $s \ge 2$; in other words, we are looking for large highly-connected multicoloured subgraphs of multicoloured graphs. When trying to work out what happens to $m(n,r,s,k)$ when $s > 1$, one quickly realises that new ideas are going to be needed. For example, none of the extremal examples from \cite{HNR1} are any help to us, since in each of them \emph{any} two colours $k$-connect almost the entire vertex set! However, we shall find that a number of the tools developed in that paper are still useful to us. We shall recall these results as we go along.

Our main results are as follows. We begin with the case $s = 2$. When also $r + 1$ is a power of $2$, we have the following fairly tight bounds.

\begin{thm}\label{r2k}
Let $n,r,k \in \N$, with $r \ge 3$, $r + 1$ a power of \textup{2} and $n \ge
16kr^2 + 4kr$. Then
$$\ds\frac{4n}{r+1} \: - \: 5kr(r + 2k + 1) \; \le \; m(n,r,2,k) \; \le \;
\ds\frac{4n}{r+1} \: + \: 4.$$
In particular, if also $k$ and $r$ are fixed, then $m(n,r,2,k) =
\ds\frac{4n}{r+1} + o(n)$.
\end{thm}

We remark that the lower bound in Theorem~\ref{r2k} in fact holds for all $r
\ge 3$, but the upper bound may increase by a factor of at most $2$ when
$r+1$ is not a power of 2. For $r = 3$ (and $n \ge 13k - 15$), however, we
can solve the problem exactly.

\begin{thm}\label{32k}
Let $n,k \in \N$, with $n \ge 13k - 15$. Then $$m(n,3,2,k) \; = \; n - k +
1.$$
\end{thm}

Our next result shows that there is a jump at $2s = r$, from $(1 - \eps)n$
to $n - 2k + 2$.

\begin{thm}\label{2ss}
Let $n,s,k \in \N$, with $n \ge \max \left\{ 2\ds{{2s} \choose s}(k-1) + 1,
13k - 15 \right\}$. Then $$m(n,2s,s,k) \; = \; n - 2k + 2.$$ Moreover, if
$2s < r \in \N$, then there exists $\eps = \eps(s,r) > 0$ such that
$$m(n,r,s,k) \; < \; (1 - \eps)n$$ for every sufficiently large $n \in
\N$.\\[-1ex]
\end{thm}

Moreover, we can determine the maximum possible $\eps$ exactly.

\begin{thm}\label{2s-1}
Let $n,s,k \in \N$ with $s \ge 2$ and $n \ge 100\ds{{2s + 1} \choose
s}^2k^2$, and let $\eps = \eps(s) = \ds{{2s+1} \choose s}^{-1}$. Then
$$(1 - \eps)n \: - \: 2{{2s + 1} \choose s}k \; \le \; m(n,2s+1,s,k) \; \le
\; (1 - \eps)n \: + \: 1.$$
\end{thm}

We have seen that (rather unsurprisingly) when $r$ is very large compared
with $s$, the function $m(n,r,s,k)/n$ is very close to $0$, and when $s$ and
$r$ are comparable the same function is close to $1$. The next theorem shows
that the function changes from one of these states to the other rather
rapidly. This is another example of a phase transition with respect to $s$.

\begin{thm}\label{jump}
Let $n,r,s,k \in \N$, with $n \ge 16kr^2 + 4kr$. Then
$$\left( 1 - e^{-s^2/4r} \right) n \: - \: 2kr\ds{r \choose \lceil s/2
\rceil} \; \le \; m(n,r,s,k) \; \le \; (s+1) \left\lceil \ds\frac{n}{
\left\lfloor \sqrt{2r} \right\rfloor} \right\rceil.$$
In particular, if $n = n(r) \gg kr\ds{r \choose \lceil s/2 \rceil}$ as $r
\to \infty$, then
$$\ds\lim_{r\rightarrow\infty}\ds\frac{m(n,r,s,k)}{n} = \left\{
\begin{array} {r@{\quad \textup{if} \quad}l} 0 & s \ll \sqrt{r} \\
1 & s \gg \sqrt{r}
\end{array}\right.$$
\end{thm}

The rest of the paper is organised as follows. In Section~\ref{tools} we
shall recall the tools developed in \cite{HNR1}; in Section~\ref{s=2} we
shall prove Theorems~\ref{r2k} and \ref{32k}; in Section~\ref{2s=r} we shall
prove Theorems~\ref{2ss} and \ref{2s-1}, and discuss the jump at $2s = r$;
in Section~\ref{s=rootr} we shall prove Theorem~\ref{jump}; and in
Section~\ref{kprobs} we shall look back on what we have learnt, and point
out some of the most obvious questions of the many that remain.

\section{Tools}\label{tools}

In \cite{HNR1} various techniques were developed to find monochromatic
$k$-connected subgraphs. Many of these tools will prove to be useful to us
below, and for the reader's convenience we shall begin by stating them here.
We start with the most crucial lemma from \cite{HNR1}, which may be proved
by induction on $m + n$.

\begin{lemma}\label{r1kbip}
Let $q,\ell,m,n \in \N$ with $m,n \ge \ell$ and $m + n \ge 2\ell + 1$. Let
$G$ be a bipartite graph with parts $M$ and $N$ of size $m$ and $n$,
respectively. If $G$ has no $(\ell+1)$-connected subgraph on at least $q$
vertices, then
\begin{eqnarray*}\label{turan}
e(G) \; \le \; \frac{q(n - \ell)(m - \ell)}{m + n - 2\ell} \: + \: (\ell^2 +
\ell)(m + n - 2\ell).
\end{eqnarray*}
\end{lemma}

We shall use Lemma~\ref{r1kbip} to prove Theorem~\ref{r2k} and
Lemma~\ref{allSbig}, below. We shall use Lemma~\ref{allSbig} to prove
Theorem~\ref{2s-1}, but we also consider the result to be interesting in its
own right. First however we note the following simple corollary of
Lemma~\ref{r1kbip}.

\begin{cor}\label{sumq}
Let $c,d,k,m,n \in \N$ with $m,n > k$, and let $G$ be a bipartite graph with
parts $M$ and $N$ of size $m$ and $n$, respectively. Let $f: E(G) \to \N$ be
a colouring of the edges of $G$, and for each $i \in \N$, let $q_k(i)$ be
the maximum order of an $k$-connected subgraph of $G$ which has all edges
coloured $i$. Suppose that all but at most $d$ edges have colours from the
set $[c]$. Then
$$\sum_{i = 1}^c q_k(i) \ge (m + n - 2k)\left(1 - \frac{d}{mn} -
\frac{ck^2(m + n)}{mn} \right).$$
\end{cor}

\begin{proof}
Apply Lemma~\ref{r1kbip} with $\ell = k - 1$ and $q = q_k(i) + 1$ for each
$i \in [c]$, and note that $$\frac{q(n - \ell)(m - \ell)}{m + n - 2\ell} \:
+ \: (\ell^2 + \ell)(m + n - 2\ell) \; \le \; \frac{qmn}{m + n - 2k} + k^2(m
+ n) - 1.$$ Adding the resulting inequalities gives $$mn - d \le \frac{mn}{m
+ n - 2k} \sum_{i = 1}^c (q_k(i) + 1) \: + \: ck^2(m + n) - c.$$ Rearranging
the inequality gives the desired result.
\end{proof}

We shall need the following observation of Bollob\'as and
Gy\'arf\'as~\cite{BG}.

\begin{obs}\label{BG}
For any graph $G$ and any $d \in \N$, either
\begin{enumerate}
\item[$(a)$] $G$ is $k$-connected, or
\item[$(b)$] $\exists \: v \in V(G)$ with $d_G(v) \le d + k - 3$, or
\item[$(c)$] $\exists \: K_{p,q} \subset \ol{G}$, with $p + q = |G| - k +
1$, and $\min\{p,q\} \ge d$.
\end{enumerate}
\end{obs}

Suppose we are given an $r$-colouring $f$ of the edges of $K_n$. For each $S
\subset [r]$, let $q_k(S)$ denote the maximum order of a $k$-connected
subgraph of $K_n$, all of the edges of which have colours from $S$.

Corollary~\ref{sumq} and Observation~\ref{BG} now allow us to prove the
following result, which we shall use in the proof of Theorem~\ref{2s-1} to
show that \emph{any} $s$-set of colours gives a large $k$-connected
subgraph. It may be thought of as a stability result for 3-colourings.

\begin{lemma}\label{allSbig}
Let $n,k,r,t \in \N$, with $n > 2t + k$, let $f$ be an $r$-colouring of
$E(K_n)$, and let $S,T,U \subset [r]$ be such that $S \cup T \cup U = [r]$.
Then either
\begin{enumerate}
\item[$(a)$] $q_k(U) \ge n - t$, or
\item[$(b)$] $q_k(S) + q_k(T) \ge n - 2t - 4k - \ds\frac{2k^2n^2}{t(n - 2t -
k)}$.
\end{enumerate}
In particular, if $q_k(U) < n - k\sqrt{n}$ and $n \ge 25k^2$, then $$q_k(S)
+ q_k(T) \ge n - 9k\sqrt{n}.$$
\end{lemma}

\begin{proof}
Let $n,k,r,t \in \N$ with $n > 2t + k$, let $f$ be an $r$-colouring of
$E(K_n)$, and let $S,T,U \subset [r]$ be such that $S \cup T \cup U = [r]$.
The result is trivial if $t \le k$, so assume that $t > k$. We divide the
problem into two cases, as follows.\\

Case 1: There exists a complete bipartite subgraph $K_{a,b}$, with $a + b
\ge n - t - k$, $b \ge a \ge t$, and all edges having colours from the set
$S \cup T$.\\

We apply Corollary~\ref{sumq} to $K_{a,b}$, with $k = k$, $c = 2$ and $d =
0$. The lemma says exactly that
\begin{eqnarray*}
q_k(S) + q_k(T) & \ge & (a + b - 2k)\left(1 - \frac{2k^2(a +
b)}{ab}\right)\\
& > & a \: + \: b \: - \: \frac{2k^2(a + b)^2}{ab} \: - \: 2k\\
& > & n -  2t  -  4k  -  \frac{2k^2n^2}{t(n - 2t - k)}
\end{eqnarray*}
and so we are done in this case.\\

Case 2: No such bipartite subgraph of $K_n$ exists, and $q_k(U) < n - t$.\\

Let $G$ be the graph with $V(G) = V(K_n)$ and $E(G) = f^{-1}(U)$, and apply
Observation~\ref{BG} to $G$ with $d = t$. Since $q_k(U) < n$, we know that
$G$ is not $k$-connected. Similarly there does not exist a complete
bipartite subgraph $K_{a,b}$ of $\ol{G}$ with $a + b = |G| - k + 1$ and $b
\ge a \ge t$, since Case 1 does not hold. Hence there must exist a vertex
$v_1 \in V(G)$ with $d_G(v_1) \le t + k - 3$.

Now let $G_1 = G - v_1 = G[V(G) \setminus \{v_1\}]$, and apply
Observation~\ref{BG} to $G_1$, again with $d = t$. Again (since  $q_k(U) < n
- 1$, and $|G_1| - k + 1 \ge n - t$), there must exist a vertex $v_2 \in
V(G_1)$ with $d_{G_1}(v_2) \le t + k - 3$. Let $G_2 = G_1 - v_2$, and
continue until we have obtained a set $X = \{v_1, \ldots, v_t\} \subset V$
satisfying $|\Gamma_G(v_i) \setminus X| \le t + k - 3$ for every $i \in
[t]$.

We now apply Corollary~\ref{sumq} to the bipartite graph with parts $X$ and
$V \setminus X$, and edges from the set $S \cup T$. The lemma says that
\begin{eqnarray*}
q_k(S) + q_k(T) & \ge & (n - 2k)\left(1 - \frac{t(t + k - 3)}{t(n-t)} -
\frac{2k^2n}{t(n-t)}\right)\\
& > & n - 2k - 2(t + k - 3) - \frac{2k^2n^2}{t(n-t)} \\
& > & n - 2t - 4k - \frac{2k^2n^2}{t(n - 2t - k)}
\end{eqnarray*}
and so we are done in this case as well.\\

The final part of the lemma follows by letting $t = \lfloor k\sqrt{n}
\rfloor$, and noting that $(k\sqrt{n} - 1)(n - 2k\sqrt{n} - k) \ge
kn^{3/2}/3$ if $n \ge 25k^2$.
\end{proof}

We shall frequently need to show that specific \emph{bipartite} graphs have
large $k$-connected subgraphs. The following observation from \cite{HNR1} is
the basic tool we use to do this.

\begin{lemma}\label{intersect}
Let $G$ be a bipartite graph with partite sets $M$ and $N$ such that $d(v)
\ge k$ for every $v \in M$, and $|\Gamma(y) \cap \Gamma(z)| \ge k$ for every
pair $y,z \in N$. Then $G$ is $k$-connected.
\end{lemma}

The next two lemmas now follow from Lemma~\ref{intersect} by removing a
suitably chosen set of `bad' vertices.

\begin{lemma}\label{21ktech}
Let $a,b,k \in \N$, and let $G$ be a bipartite graph with parts $M$ and $N$
such that $|M| \ge 4b + k$, $|N| \ge a \ge 2k$, and $d(v) \ge |M| - b$ for
every $v \in N$. Then $G$ has a $k$-connected subgraph on at least $$|G| -
\frac{ab}{a - k + 1} \; > \; |G| - 2b$$ vertices.
\end{lemma}

\begin{proof}
Let $a,b \in \N$ and $G$ be as described. Let $$U = \{v \in M : d_G(v) \le k
- 1\},$$ and observe that each vertex of $U$ sends at least $|N| - k + 1$
non-edges into $N$, and that there are in total at most $b|N|$ non-edges
between $M$ and $N$ (since $d(v) \ge |M| - b$ for every $v \in N$). Hence
$$|U|(|N| - k + 1) \; \le \; b|N|,$$ and so $$|U| \; \le \;
\ds\frac{b|N|}{|N| - k + 1} \; \le \; \frac{ab}{a - k + 1} \; < \; 2b,$$
since the function $\ds\frac{bx}{x - k + 1}$ is decreasing for $x > k - 1$,
and $|N| \ge a \ge 2k$. Now, consider the bipartite graph $G' = G[M
\setminus U, N]$. Each vertex of $M \setminus U$ has degree at least $k$ in
$G'$, by the definition of $U$, and each pair of vertices of $N$ have at
least $k$ common neighbours in $M$, since $|M| \ge 4b + k$, so $|M \setminus
U| \ge 2b + k$, and each vertex of $N$ has at most $b$ non-neighbours in
$M$. Thus, by Lemma~\ref{intersect}, $G'$ is $k$-connected, and has order
$$|G| - |U| \; \ge \; |G| - \frac{ab}{a - k + 1} \; > \; |G| - 2b.$$
\end{proof}

The following easy lemma is very similar to Lemma 15 of \cite{HNR1}, but a
little stronger. In particular, we have removed the requirement that $3|M|
\ge |N|$.

\begin{lemma}\label{2s-1bip}
Let $k \in \N$, and let $G$ be the complete bipartite graph with parts $M$
and $N$, where $|N| \ge |M| \ge 15k$. Let $f$ be an $r$-colouring of $E(G)$,
and let $S,T,U \subset [r]$ be such that $S \cup T \cup U = [r]$. Suppose
that
\begin{enumerate}
\item[$(a)$] $|\{v \in N \: : f(uv) \in S\}| \le k$ for every $u \in M$, and
\item[$(b)$] $|\{u \in M : f(uv) \in T\}| \le k$ for every $v \in N$.
\end{enumerate} Then there exists a $k$-connected subgraph of $G$, using
only colours from $U$, and avoiding at most $5k$ vertices of $M$ and $2k$
vertices of $N$. In particular, $q_k(U) \ge |G| - 7k$.
\end{lemma}

\begin{proof}
We may assume that $r = 3$, and that $S = \{1\}$, $T = \{2\}$ and $U =
\{3\}$. Let $k,m,n \in \N$ with $n \ge m \ge 15k$, let $|M| = m$ and $|N| =
n$, and let $f$ be a 3-colouring of $E(G)$ satisfying the conditions of the
lemma. Let
$$S_M = \{v \in M : v\textup{ sends at most }3n/5\textup{ edges of colour
}3\textup{ into }N\}\textup{, and}$$
$$S_N = \{v \in N : v \textup{ sends at most } 6k \textup{ edges of colour
}3\textup{ into }M\} \hspace{1.3cm}$$
be sets of `bad' vertices. We shall remove the bad sets and apply
Lemma~\ref{intersect}.

We need to bound $|S_M|$ and $|S_N|$ from above. Since each vertex of $M$
has at most $k$ incident edges of colour $1$, we have $|f^{-1}(1)| \le km$,
and similarly $|f^{-1}(2)| \le kn$. Also, since each vertex of $S_M$ has at
least $2n/5$ incident edges of colour $1$ or $2$, we have $|f^{-1}(1)| +
|f^{-1}(2)| \ge |S_M|(2n/5)$. Finally, each vertex of $S_N$ has at most $6k$
incident edges of colour 3, and at most $k$ incident edges of colour 2, so
has at least $(m - 7k)$ incident edges of colour 1. Hence $|f^{-1}(1)| \ge
|S_N|(m - 7k)$. Thus
$$|S_M| \; \le \; \frac{5}{2n}\Big( |f^{-1}(1)| + |f^{-1}(2)| \Big) \; \le
\; \frac{5k(m+n)}{2n} \; \le \; 5k, \: \textup{ and}$$
$$|S_N| \; \le \; \frac{|f^{-1}(1)|}{m - 7k} \; \le \; \frac{km}{m - 7k} \;
\le \; 2k, \hspace{0.8cm}$$ since $m \ge 14k$.

Now, let $M' = M \setminus S_M$ and $N' = N \setminus S_N$, and let $H$ be
the bipartite graph with vertex set $M' \cup N'$, and edge set $f^{-1}(3)$.
If $x \in N'$, then $x$ sends at least $6k$ edges of colour $3$ into $M$, so
$$d_H(x) \; \ge \; 6k - |S_M| \; \ge \; 6k - 5k \; = \; k,$$ and similarly
if $y,z \in M'$, then
\begin{eqnarray*}
|\Gamma_H(y) \cap \Gamma_H(z)| & \ge & 3n/5 \: + \: 3n/5 \: - \: n
\: - \: |S_N|\\[+1ex]
& = & n/5 - |S_N| \; \ge \; 3k - 2k \; = \; k,
\end{eqnarray*}
since $n \ge 15k$, so the conditions of Lemma~\ref{intersect} are satisfied.
Thus by Lemma~\ref{intersect}, $H$ is $k$-connected. Since also $|M
\setminus V(H)| = |S_M| \le 5k$ and $|N \setminus V(H)| = |S_N| \le 2k$, $H$
is the desired subgraph.
\end{proof}

The results above will be our main tools in the sections that follow.
However, we shall also use the following well-known theorem of
Mader~\cite{Mader} in the proofs of Theorems~\ref{r2k} and \ref{jump}.

\begin{mader}
Let $\alpha \in \RR$, and let $G$ be a graph with average degree $\alpha$.
Then $G$ has an $\alpha/4$-connected subgraph.
\end{mader}

Mader's Theorem implies that a monochromatic $(n-1)/4r$-connected subgraph
exists in any $r$-colouring of $E(K_n)$ (to see this, simply consider the
colour which is used most frequently). This subgraph is $k$-connected if $n
\ge 4kr + 1$, and has at least $(n-1)/4r + 1$ vertices. It is this weak
bound that we shall use.

We also state the following result from \cite{HNR1} here, so that we may
refer to it more easily. We shall use Theorem~\ref{21k} to prove the lower
bounds in Theorems~\ref{32k} and \ref{2ss}.

\begin{thm}{\label{21k}}
Let $n,k \in \N$, with $n \ge 13k - 15$. Then $$m(n,2,1,k) \; = \; n - 2k +
2.$$
\end{thm}

Finally we make some simple observations about $k$-connected graphs.

\begin{obs}\label{addvtx}
Let $G$ be a graph, and $v \in V(G)$. If $G - v$ is $k$-connected and $d(v)
\ge k$, then $G$ is also $k$-connected.
\end{obs}

Given graphs $G$ and $H$, define $G \cup H$ to be the graph with vertex set
$V(G) \cup V(H)$ and edge set $E(G) \cup E(H)$.

\begin{obs}\label{linkup}
Let $k \in \N$, and $H_1$ and $H_2$ be $k$-connected subgraphs of a graph
$G$. If there exist $k$ vertices $\{v_1, \ldots, v_k\} \subset V(H_1)$ such
that $|\Gamma(v_i) \cap V(H_2)| \ge k$ for each $i \in [k]$, then $H_1 \cup
H_2$ is $k$-connected.
\end{obs}

\begin{obs}\label{GcupH}
Let $k \in \N$. If $G$ and $H$ are $k$-connected graphs, and $|V(G) \cap
V(H)| \ge k$, then the graph $G \cup H$ is also $k$-connected.
\end{obs}

Throughout, we shall write $V$ for $V(K_n)$. For any undefined terms see
either \cite{MGT} or \cite{HNR1}.\\

\section{The case $s = 2$}\label{s=2}

We begin at the bottom, with the case $s = 2$. We shall be able to give
fairly tight bounds on $m(n,r,2,k)$ for infinitely many value of $r$. We
begin with a construction, which will give us our upper bound.

\begin{lemma}\label{henry}
For every $n,r,s,k \in \N$, we have $$m(n,r,s,k) \: \le \: m(n,r,s,1) \: \le
\: 2^s \left\lceil \ds\frac{n}{2^{\lfloor \log_2(r+1) \rfloor}} \right\rceil
\: < \: 2^s \left\lceil \ds\frac{2n}{r+1} \right\rceil.$$ In particular, if
$r + 1$ is a power of \textup{2} and a divisor of $n$, then $$m(n,r,s,k) \le
\ds\frac{2^sn}{r + 1}.$$
\end{lemma}

\begin{proof}
Let $n,r,s,k \in \N$, and let $R = \lfloor \log_2(r + 1) \rfloor$. We shall
define a $(2^R - 1)$-colouring $f: E(K_n) \to \{0,1\}^R \setminus
\textbf{0}$ of the edges of $K_n$. First partition $V = V(K_n)$ into $2^R$
subsets $\{V_x : x \in \{0,1\}^R\}$ of near-equal size (i.e., $|(|V_x| -
|V_y|)| \le 1$ for every pair $x, y$). Now if $i \in V_x$, $j \in V_y$ and
$x \neq y$, then let $f(ij) = x - y \pmod 2$; for the remaining edges choose
$f$ arbitrarily.

Choose a subset $S \subset \{0,1\}^R \setminus \textbf{0}$ of size $s$, a
vector $x \in \{0,1\}^R$ and a vertex $v \in V_x$. Let $G$ be the graph with
vertex set $V$ and edge set $f^{-1}(S)$, and let $\P(S) = \{y_1 + \ldots +
y_t \pmod 2 \in \{0,1\}^R : t \in \N$, and $y_i \in S$ for each $i \in
[t]\}$. Note that $|\P(S)| \le 2^s$.

Now, there is a path from $v$ to a vertex $u \in V_y$ using only edges of
$S$ if and only if $x - y \in \P(S)$, since such a path corresponds to a sum
of vectors from $S$. Hence the component of $v$ in $G$ is exactly
$\bigcup\{V_y : x - y \in \P(S)\}$.

Since $v$ and $S$ were arbitrary, $|V_y| \le \lceil n/2^R \rceil$ for each
$y \in \{0,1\}^R$, and $|\P(S)| \le 2^s$, it follows that in the colouring
$f$, there is no $1$-connected subgraph using at most $s$ colours on more
than $$2^s \left\lceil \ds\frac{n}{2^R} \right\rceil = 2^s \left\lceil
\ds\frac{n}{2^{\lfloor \log_2(r+1) \rfloor}} \right\rceil$$ vertices. This
proves that $m(n,r,s,1) \le 2^s \left\lceil \ds\frac{n}{2^{\lfloor
\log_2(r+1) \rfloor}} \right\rceil$; the remaining inequalities are trivial.
\end{proof}

The following corollary is immediate from the lemma; we state it just for
emphasis.

\begin{cor}\label{r2kupp}
If $n,r,k \in \N$, and $r + 1$ is a power of \textup{2}, then $$m(n,r,2,k)
\: \le \: 4 \left\lceil \ds\frac{n}{r+1} \right\rceil.$$
\end{cor}

We shall prove an almost matching lower bound on $m(n,r,2,k)$. To help the
reader (and because we can prove a stronger result in this case), we begin
with the case $k = 1$. We shall need the following result from \cite{HNR1}.

\begin{lemma}\label{r11bip}
The order of the largest monochromatic component of an $r$-colouring of
$E(K_{m,n})$ is at least $\ds\frac{m+n}{r}$.
\end{lemma}

\begin{thm}\label{r21}
Let $n,r \in \N$, with $r \ge 3$. Then $$m(n,r,2,1) \ge \ds\frac{4n}{r+1}.$$
\end{thm}

\begin{proof}
Let $n,r \in \N$, with $r \ge 3$, and let $f$ be an $r$-colouring of
$E(K_n)$. We shall show that there exists a connected subgraph of $K_n$, on
at least $\ds\frac{4n}{r+1}$ vertices, using only at most two colours in the
colouring $f$. For aesthetic reasons, we shall assume that $r + 1 \: | \:
4n$ (otherwise the proof is almost the same, but slightly messier). Let $G$
be the largest connected \emph{monochromatic} subgraph of $K_n$, let $A =
V(G)$, and let $|A| = \ds\frac{cn}{r+1}$. If $c \ge 4$ then we are done, so
assume that $c < 4$.

Let $B = V \setminus A$, and without loss of generality, assume that the
edges of $G$ all have colour $1$. Then since $G$ is maximal, no edge in the
bipartite graph $H = K_n[A,B]$ has colour $1$, so some colour occurs at
least $\ds\frac{|A||B|}{r-1}$ times in $H$. Again without loss let this
colour be $2$, and let $B_2 = \{u \in B : f(uv) = 2$ for some $v \in A\}$ be
the set of vertices of $B$ which are incident to some edge of $H$ of colour
$2$.

Suppose first that $|B_2| \ge \ds\frac{(4 - c)n}{r + 1}$. Then the set $A
\cup B_2$ is connected by colours $1$ and $2$, and $|A \cup B_2| \ge
\ds\frac{4n}{r + 1}$, so we are done. So assume that $|B_2| < \ds\frac{(4 -
c)n}{r + 1}$, and choose a set $B_2'$ such that $B_2 \subset B_2' \subset
B$, and $|B_2'| = \ds\frac{(4 - c)n}{r + 1}$.

We apply Lemma~\ref{r11bip} to the bipartite graph $H_2$ with parts $A$ and
$B_2'$, and edges of colour 2. We have
\begin{eqnarray*}
e(H_2) & \ge & \ds\frac{|A|\,|B|}{r-1} \; = \; \Big( \ds\frac{1}{r-1} \Big)
\Big( \ds\frac{cn}{r+1} \Big) \Big( n - \ds\frac{cn}{r+1} \Big)\\[+1ex] & =
& \Big( \ds\frac{r+1-c}{(4-c)(r-1)} \Big) \Big( \ds\frac{cn}{r+1} \Big)
\Big( \ds\frac{(4-c)n}{r+1} \Big)\\[+1ex]
& = & \Big( \ds\frac{r+1-c}{(4-c)(r-1)} \Big) \: |A| \: |B_2'|,
\end{eqnarray*} so Lemma~\ref{r11bip} implies that there exists a connected
subgraph of $H_2$ on at least
\begin{eqnarray*}
\ds\frac{e(H_2)(|A| + |B_2'|)}{|A|\,|B_2'|} \; \ge \; \Big( \ds\frac{(r + 1
- c)}{(4-c)(r-1)} \Big) \Big( \ds\frac{4n}{r+1} \Big)
\end{eqnarray*}
vertices, since $|A| + |B_2'| = \ds\frac{4n}{r+1}$. This subgraph is
monochromatic, and so, since $G$ was chosen to be the largest monochromatic
subgraph, we have
$$\ds\frac{4(r + 1 - c)n}{(4-c)(r-1)(r+1)} \; \le \; \ds\frac{cn}{r+1},$$
which implies that
$$(r-1)c^2 - 4rc + 4(r + 1) \le 0,$$ since $c < 4$ and $r > 1$. The
quadratic factorises as $(c - 2)((r-1)c - 2(r + 1)) \le 0$, so we have $2
\le c \le \ds\frac{2(r+1)}{r-1}$.

We shall only need that $c \ge 2$. For suppose some vertex $u \in B$ sends
edges of only one colour into $A$, i.e., $|\{i \in [r] : f(uv) = i$ for some
$v \in A\}| = 1$. Let that colour be $j$, and consider the star, centred at
$u$, with edges of colour $j$. It is monochromatic, connected, and has order
larger than $G$, a contradiction. Thus every vertex in $B$ sends edges of at
least two different colours into $A$, and so, by the pigeonhole principle,
some colour ($\ell$, say) is sent by at least $\ds\frac{2|B|}{r-1}$
different vertices of $B$.

Let $D = \{u \in B : f(uv) = \ell$ for some $v \in A\}$. We have
\begin{eqnarray*}
|A \cup D| & \ge & \ds\frac{cn}{r+1} + \ds\frac{2|B|}{r - 1} \; = \;
\ds\frac{n}{r + 1} \left( c + \ds\frac{2(r + 1 - c)}{r - 1} \right)\\[+1ex]
& = & \ds\frac{n}{r+1} \left( \ds\frac{c(r - 3) + 2(r + 1)}{r - 1}
\right)\\[+1ex]
& \ge & \ds\frac{n}{r+1} \left( \ds\frac{2(r - 3) + 2(r + 1)}{r - 1} \right)
\; = \; \ds\frac{4n}{r+1},
\end{eqnarray*}
the last inequality following because $c \ge 2$ and $r \ge 3$. The vertices
of $A \cup D$ are connected by edges
of colour $1$ and $\ell$, so we are done.
\end{proof}

We now modify the proof of Theorem~\ref{r21} to prove Theorem~\ref{r2k}. We
shall use Mader's Theorem, and Lemmas~\ref{r1kbip} and \ref{21ktech}.

\begin{proof}[Proof of Theorem~\ref{r2k}]
Let $n,r \in \N$, with $r \ge 3$ and $n \ge 16kr^2 + 4kr$. The upper bound
follows by Lemma~\ref{henry} if $r + 1$ is a power of 2; it remains to prove
the lower bound. If $n < kr(r+1)(r + 2k + 1)$ then the result is trivial, so
assume $n \ge kr(r+1)(r + 2k + 1)$.

Let $f$ be an $r$-colouring of the edges of $K_n$. We shall show that there
exists a connected subgraph of $K_n$, on at least $\ds\frac{4n}{r+1} - 5kr(r
+ 2k + 1)$ vertices, using only at most two colours in the colouring $f$.
Let $G$ be the largest $k$-connected \emph{monochromatic} subgraph of $K_n$,
let $A = V(G)$, and let $|A| = \ds\frac{cn}{r+1}$. Following the proof of
Theorem~\ref{r21}, we shall show that $c \ge 2 - \ds\frac{17kr^2(r +
2k)}{n}$. Suppose for a contradiction that $c < 2 - \ds\frac{17kr^2(r +
2k)}{n}$.

Let $B = V \setminus A$, and without loss of generality, assume that the
edges of $G$ all have colour $1$. Now, since $G$ is maximal, no vertex in
$B$ sends more than $k - 1$ edges of colour $1$ into $A$, by
Observation~\ref{addvtx}. Thus in the bipartite graph $H = K_n[A,B]$, some
colour occurs at least $\ds\frac{(|A| - k + 1)|B|}{r-1}$ times; without loss
let this colour be $2$. Let $B_2 = \{u \in B : |\{v \in A : f(uv) = 2\}| \ge
k \}$ be the set of vertices of $B$ which are incident to at least $k$ edges
of $H$ of colour $2$.

Suppose first that $|B_2| \ge \ds\frac{(4 - c)n}{r + 1}$. Then the set $A
\cup B_2$ is $k$-connected by colours $1$ and $2$ by
Observation~\ref{addvtx}, and $|A \cup B_2| \ge \ds\frac{4n}{r + 1}$, so we
are done. So assume that $|B_2| < \ds\frac{(4 - c)n}{r + 1}$, and choose a
set $B_2'$ such that $B_2 \subset B_2' \subset B$, and $|B_2'| =
\left\lfloor \ds\frac{(4 - c)n}{r+1} \right\rfloor$.

We shall apply Lemma~\ref{r1kbip} to the bipartite graph $H_2$ with parts
$A$ and $B_2'$, and edges of colour 2. First note that
\begin{equation}\label{eH2} e(H_2) \; \ge \; \ds\frac{(|A| -
r(k-1))|B|}{r-1}, \end{equation}
since we discarded at most $(k-1)|B|$ edges of colour 2 from $H$ when
forming $H_2$. Let $\ell = k - 1$; we must check that $|A|, |B_2'| \ge \ell$
and $|A| + |B_2'| \ge 2\ell + 1$. These bounds follow because $n > 4kr$, so
$$|A| \; \ge \; \ds\frac{n - 1}{4r} + 1 \; > \; k$$
by Mader's Theorem, and
$$|B_2'| \; \ge \; \ds\frac{(4 - c)n}{r+1} - 1 \; > \; \frac{2n}{r-1} - 1 \;
 > \; k,$$
since we assumed that $c < 2$.

So, by Lemma~\ref{r1kbip}, if there does not exist a $k$-connected subgraph
of $H_2$ on at least $q$ vertices, then
\begin{eqnarray}
e(H_2) & \le & \ds\frac{q(|A| - \ell)(|B_2'| - \ell)}{|A| + |B_2'| - 2\ell}
\: + \: (\ell^2 + \ell)(|A| + |B_2'| - 2\ell)\nonumber\\[+1ex]
& \le & \ds\frac{ q \Big( \ds\frac{cn}{r+1} \Big) \Big( \ds\frac{(4-c)n}
{r+1} \Big) }{ \ds\frac{4n}{r+1} - 2k } \: + \: (\ell^2 + \ell) \left(
\ds\frac{4n}{r+1} \right)\nonumber\\[+1ex]
& \le & \ds\frac{qc(4-c)n^2} {(r+1)(4n - 2k(r+1))} \: + \: k^2 \left(
\ds\frac{4n}{r+1}\right),\label{q}
\end{eqnarray}
since $|A| + |B_2'| \ge \ds\frac{4n}{r+1} - 1$ and $k = \ell + 1$. Combining
\eqref{eH2} and \eqref{q}, we get
\begin{eqnarray}
q & \ge & \ds\frac{(r+1)(4n - 2k(r+1))}{c(4-c)n^2} \left( \ds\frac{\big(|A|
- r(k-1)\big)|B|}{r-1} - k^2\left( \ds\frac{4n}{r+1} \right)
\right)\nonumber \\[+1ex]
& \ge & \ds\frac{4n(r+1) - 2k(r+1)^2}{c(4-c)n^2} \left( \ds\frac{\big(cn -
kr(r+1)\big)(r+1-c)n}{(r-1)(r+1)^2} - \ds\frac{4k^2n}{r+1}
\right)\nonumber\\[+1ex]
& = & \ds\frac{4(r+1-c)n}{(4-c)(r-1)(r+1)} \: - \: \left(
\ds\frac{4kr(r+1-c)}{c(4-c)(r-1)} + \ds\frac{16k^2}{c(4-c)} +
\ds\frac{2k(r+1-c)}{(4-c)(r-1)} \right) \nonumber \\[+1ex]
&    & \hspace{3.9cm} + \; \ds\frac{1}{n} \left( \ds\frac{ 2k^2r(r+1)(r+1-c)
}{ c(4-c)(r-1) } + \frac{8k^3(r+1)}{c(4-c)} \right)\nonumber \\[+1ex]
& > & \ds\frac{4(r+1-c)n}{(4-c)(r-1)(r+1)} \: - \: 17k(r + 2k),\label{q2}
\end{eqnarray}
since $\ds\frac{r+1-c}{r-1} < 2$ and $\ds\frac{1}{4} < c < 2$, so $c(4-c) >
\ds\frac{1}{2}$.

Inequality \eqref{q2} holds if there does not exist a $k$-connected subgraph
of $H_2$ on at least $q$ vertices. Therefore, since $G$ was chosen to be the
largest $k$-connected monochromatic subgraph of $K_n$, and $|G| =
\ds\frac{cn}{r+1}$, we have
$$\ds\frac{cn}{r+1} \: + \: 1 \; > \; \ds\frac{4(r+1-c)n} {(4-c)(r-1)(r+1)}
\: - \: 17k(r + 2k),$$ which implies
$$(c - 2)((r-1)c - 2(r + 1)) \; < \; \ds\frac{68kr^2(r + 2k)}{n},$$ as in
the proof of Theorem~\ref{r21}. Now, $c < 2$, so $(r-1)c - 2(r+1) < -4$, and
thus
$$2 - c \; < \; \ds\frac{17kr^2(r + 2k)}{n},$$ so $c \ge 2 -
\ds\frac{17kr^2(r + 2k)}{n}$, as claimed.

Now, for each $j \in [r]$ define $$C_j = \{v \in B : |\{u \in A : f(uv) \neq
j\}| \le kr\},$$ and suppose that $|C_j| \ge 2kr$ for some colour $j \in
[r]$, i.e., there are at least $2kr$ distinct vertices in $B$ which each
send at least $|A| - kr$ edges of colour $j$ into $A$.

We shall apply Lemma~\ref{21ktech} to obtain a contradiction. Let $F$ be the
bipartite graph with parts $A$ and $C_j$, and edges of colour $j$, and let
$a = 2kr$ and $b = kr$. Now, $|C_j| \ge 2kr \ge 2k$, and $|A| \ge \ds\frac{n
- 1}{4r} + 1 \ge 4kr + k$ (by Mader's Theorem, and because $n \ge 16kr^2 +
4kr$), and $d_F(v) \ge |A| - kr$ for every $v \in C_j$, by the definition of
$C_j$.

Thus by Lemma~\ref{21ktech}, there exists a $k$-connected subgraph $F'$ of
$F$ on more than $|A| + |C_j| - 2b$ vertices. But $|C_j| \ge 2kr = 2b$, so
$|V(F')| > |A|$, and $F'$ is monochromatic. This is a contradiction, since
$G$ was chosen to be the largest monochromatic $k$-connected subgraph of
$K_n$.

So for each colour $j \in [r]$, there are at most $2kr$ distinct vertices in
$B$ which send at least $|A| - kr$ edges of colour $j$ into $A$. We remove
these vertices from $B$ to obtain $$B' = \{v \in B : |\{u \in A : f(uv) \neq
j\}| > kr\textup{ for every }j \in [r]\},$$ with $|B'| \ge |B| - 2kr(r-1)$
(note that $|C_1| = 0$). Now, for each vertex $v \in B'$, we have $|\{i \in
[r] : |\{u \in A : f(uv) = i\}| \ge k\}| \ge 2$, i.e., $v$ sends at least
$k$ edges of at least two different colours into $A$. Therefore, by the
pigeonhole principle, there must exist a colour, $\ell$ say, such that at
least $\ds\frac{2|B'|}{r-1}$ vertices of $B'$ send at least $k$ edges of
colour $\ell$ into $A$.

Let $D = \{v \in B : |\{u \in A : f(uv) = \ell\}| \ge k\}$. We
have
\begin{eqnarray*}
|A \cup D| & \ge & \ds\frac{cn}{r+1} \: + \: \ds\frac{2|B'|}{r - 1} \; \ge
\; \ds\frac{cn}{r+1} \: + \: \ds\frac{2|B|}{r - 1} \: - \: 4kr\\[+1ex]
& = & \ds\frac{n}{r + 1} \left( c + \ds\frac{2(r + 1 - c)}{r - 1} \right) \:
- \: 4kr\\[+1ex]
& = & \ds\frac{n}{r+1} \left( \ds\frac{c(r - 3) + 2(r + 1)}{r - 1} \right)
\: - \: 4kr\\[+1ex]
& \ge & \ds\frac{4n}{r+1} \: - \: \ds\frac{17kr^2(r + 2k)(r - 3)}{4(r - 1)(r
+ 1)} \: - \: 4kr\\[+1ex]
& > & \ds\frac{4n}{r+1} \: - \: 5kr(r + 2k + 1)
\end{eqnarray*} since $c \ge 2 - \ds\frac{17kr^2(r + 2k)}{n}$ and $r \ge 3$.
Now by Observation~\ref{addvtx}, the subgraph of $K_n$ with vertex
set $A \cup D$ and all edges of colour $1$ or $\ell$ is
$k$-connected, so we are done.
\end{proof}

When $r = 3$ we can do better than Theorem~\ref{r2k}; in fact we can
determine the function $m(n,3,2,k)$ exactly when $n \ge 13k - 15$. The alert
reader will have noticed that this is the same bound on $n$ as we obtained
in Theorem~\ref{21k} -- this is not coincidence, the bound is necessary
because we shall use Theorem~\ref{21k} in the proof of Theorem~\ref{32k}!

The following simple construction gives us our upper bound.

\begin{lemma}\label{32kupp}
Let $n,k \in \N$. If $n \le 3k - 3$ then $m(n,3,2,k) = 0$. If $n \ge 3k-2$
then $m(n,3,2,k) \le n - k + 1$.
\end{lemma}

\begin{proof}
Let $n,k \in \N$ with $n \ge 3k - 2$. Let $A$, $B$ and $C$ be pairwise
disjoint subsets of $V = V(K_n)$, each of size $k-1$, and let $W = V
\setminus (A \cup B \cup C)$. Colour the edges between $A$ and $B \cup W$
with colour 1, those between $B$ and $C \cup W$ with colour 2, and those
between $C$ and $A \cup W$ with colour 3. Colour the edges inside the sets
arbitrarily.

Let $H$ be a $k$-connected subgraph on at least $n - k + 1$ vertices, using
at most two colours, and let these colours be $1$ and $2$ (the proof in the
other cases is identical). Let $V(H) = X$. Since $n \ge 3k-2$, $|X| \ge 2k -
1$, so the set $X \cap (A \cup W)$ is non-empty. Let $u \in X \cap (A \cup
W)$. Now $X \cap C = \emptyset$, since if $v \in X \cap C$, then $u$ and $v$
are disconnected in $H[X \setminus B]$, which is a contradiction, since $|X
\cap B| \le |B| = k - 1$. Since $X \cap C = \emptyset$ and $|C| = k - 1$, we
have $|X| \le n - k + 1$.

We have shown that in the colouring described above, there is no
$k$-connected subgraph using only two colours on more than $n - k + 1$
vertices. Therefore $m(n,3,2,k) \le n - k + 1$ when $n \ge 3k - 2$.

Now let $n,k \in \N$ with $n \le 3k - 3$. Partition $V$ into parts $A$, $B$
and $C$, each of size at most $k-1$, and colour the edges between the parts
as above: colour 1 between $A$ and $B$, colour 2 between $B$ and $C$, and
colour 3 between $C$ and $A$. This time, however, colour the edges inside
$A$ with colour 2, those inside $B$ with colour 3, and those inside $C$ with
colour 1. Now it is easy to check that there is no $k$-connected subgraph
using only two colours, so $m(n,3,2,k) = 0$ as claimed.
\end{proof}

We now prove the matching lower bound when $n \ge 13k - 15$. The argument is
similar to the proof of Theorem~\ref{21k} in \cite{HNR1} -- just one extra
idea is needed.

\begin{proof}[Proof of Theorem~\ref{32k}]
The upper bound follows from Lemma~\ref{32kupp}, and for $k = 1$ the result
is trivial, so let $n,k \in \N$ with $k \ge 2$ and $n \ge 13k - 15$, and let
$f$ be a $3$-colouring of the edges of $K_n$. We shall find a $k$-connected
subgraph $H$ of $K_n$, using at most $2$ colours of $f$, on at least $n - k
+ 1$ vertices.

For $i = 1,2,3$, let $G^{(i)}$ denote the graph with vertex set $V = V(K_n)$
and edge set $f^{-1}(i)$ (the edges of colour $i$), and for each pair
$\{i,j\} \subset \{1,2,3\}$, let $G^{(i,j)}$ denote the subgraph with vertex
set $V$ and edge set $f^{-1}(i) \cup f^{-1}(j)$ (the edges of colour $i$ or
$j$).

We shall first find two $k$-connected subgraphs, using at most two colours
each, which cover the vertex set $V$. Since $n \ge 13k - 15$, by
Theorem~\ref{21k} either $G^{(1,2)}$ or $G^{(3)}$ contains a $k$-connected
subgraph $H$ on at least $n - 2k + 2 \ge 11k - 13 > 2k - 1$ vertices.
Suppose that $H$ is in $G^{(3)}$, and let $V(H) = X$. Let $A$ be the set of
vertices of $V \setminus X$ which send at least $k$ edges of colour 1 or 3
into $X$, and let $B$ be the set of vertices of $V \setminus X$ which send
at least $k$ edges of colour 2 or 3 into $X$. Since $|X| \ge 2k - 1$, we
have $A \cup B = V \setminus X$. Without loss of generality, let $|A| \ge
|B|$. Now $G^{(1,3)}[X \cup A]$ is $k$-connected, by
Observation~\ref{addvtx}, and
$$|X \cup A| \; \ge \; |X| + \ds\frac{n - |X|}{2} \; \ge \; n - k + 1,$$
since $|X| \ge n - 2k + 2$, so we have found the desired subgraph.

So we may assume that $G^{(1,2)}$ contains a $k$-connected subgraph on at
least $n - 2k + 2$ vertices, and similarly for $G^{(1,3)}$ and $G^{(2,3)}$.
Let $Y$ be the vertex set of the largest $k$-connected subgraph in
$G^{(1,2)}$, and let $Z$ be the vertex set of the largest $k$-connected
subgraph in $G^{(1,3)}$. Since $|Y|,|Z| \ge n - 2k + 2$, $n \ge 13k - 15$
and $k \ge 2$ we have
$$|Y \cap Z| \; \ge \; n - 4k + 4 \; \ge \; 9k - 11 \; > \; 2k - 1.$$

We claim that $Y \cup Z = V$. To see this, suppose there is a vertex $v \in
V \setminus (Y \cup Z)$. Since $|Y \cap Z| \ge 2k - 1$, $v$ must send at
least $k$ edges of colour 1 or 2, or at least $k$ edges of colour 1 or 3
into $Y \cap Z$. Without loss of generality, assume that $v$ sends at least
$k$ edges of colour 1 or 2. Then $G^{(1,2)}[Y \cup \{v\}]$ is $k$-connected
by Observation~\ref{addvtx}, contradicting the maximality of $Y$. So $Y \cup
Z = V$, as claimed, and we have found two $k$-connected bichromatic
subgraphs which cover $V$.

Now, let $C = Y \setminus Z$, and $D = Z \setminus Y$. If $|Z| \ge n - k +
1$ then $G^{(1,3)}[Z]$ is the desired $k$-connected subgraph, so assume not.
Therefore $|C| \ge k$, and similarly we may assume that $|D| \ge k$. We wish
to apply Lemma~\ref{intersect} to the bipartite graph $G' = G^{(2,3)}[Y \cap
Z, C \cup D]$, so let $M' = Y \cap Z$ and $N = C \cup D$. We must first
remove the `bad' vertices, of degree at most $k-1$ in $G'$, from the graph.
As in the proof of Lemma~\ref{21ktech}, define $$U = \{v \in M' : d_{G'}(v)
\le k - 1\}.$$ We shall show that $|U| \le k - 1$.

For each $i \in \{1,2,3\}$, let $r(i) = |f^{-1}(i) \cap E(C,D)|$ be the
number of edges between $C$ and $D$ that are coloured $i$. Since $Z$ is
maximal, each vertex of $C$ can send at most $k-1$ edges of colour 1 or 3
into $Z$, so $G^{(1,3)}[C,Z]$ has at most $|C|(k - 1)$ edges. Therefore,
$G^{(1)}[C,Y \cap Z]$ has at most $|C|(k - 1) - r(1) - r(3)$ edges.
Similarly, $G^{(1)}[D,Y \cap Z]$ has at most $|D|(k - 1) - r(1) - r(2)$
edges, so $G'$ has at most $$|N|(k - 1) - |C||D| - r(1)$$ \emph{non}-edges,
since $|C| + |D| = |N|$ and $\sum_{i=1}^3 r(i) = |C||D|$.

Now, by the definition of $U$, each vertex of $U$ sends at least $|N| - k +
1$ edges of colour 1 into $N = C \cup D$, so $G'$ has at least $|U|(|N| - k
+ 1)$ non-edges. Hence
\begin{eqnarray*}
|U|(|N| - k + 1) & \le & |N|(k - 1) - |C||D| - r(1)\\[+1ex] & \le & |N|(k -
1) - k^2,
\end{eqnarray*}
since $|C|, |D| \ge k$ and $r(1) \ge 0$. Thus
$$|U| \; \le \; \ds\frac{|N|(k-1) - k^2}{|N|-k+1} \; = \; k - 1 -
\ds\frac{2k - 1}{|N| - k + 1} \; < \; k - 1,$$
since $|N| - k + 1 > 0$.

We complete the proof of Theorem~\ref{32k} by setting $M = M' \setminus U$,
and applying Lemma~\ref{intersect} to the graph $G = G^{(2,3)}[M,N]$. By the
definition of $U$, $d_G(x) \ge k$ for every vertex $x \in M$, and
\begin{eqnarray*}
|M| & = & |Y \cap Z| - |U| \; \ge \; (n - 4k + 4) - (k - 1)\\[+1ex]
& = & n - 5k + 5 \; \ge \; 8k - 10 \; \ge \; 3k - 2,
\end{eqnarray*}
since $|U| \le k - 1$, $n \ge 13k - 15$ and $k \ge 2$. Also $d_G(y) \ge |M|
- k + 1$ for every $y \in N$, since each vertex of $C \cup D$ sends at most
$k - 1$ edges of colour 1 to $Y \cap Z$. Therefore,
$$|\Gamma_G(y) \cap \Gamma_G(z)| \; \ge \; |M| - 2k + 2 \; \ge \; k$$
for every pair $y,z \in N$, so by Lemma~\ref{intersect}, $G$ is
$k$-connected.

Since $M \cup N = V \setminus U$ and $|U| \le k - 1$, $G$ is the desired
$k$-connected subgraph using at most two colours.
\end{proof}

\begin{rmk}
We needed the bound $n \ge 13k - 15$ in order to apply Theorem~\ref{21k} --
the rest of the proof required only $n \ge 8k - 7$. Therefore any
improvement on the bound on $n$ in Theorem~\ref{21k} would give an immediate
improvement here also.
\end{rmk}

When $r + 1$ is not a power of 2, we can in general only determine
$m(n,r,2,k)$ up to a factor of 2. However, it follows from
Theorem~\ref{2s-1}, which we shall prove in the next section, that
$m(n,5,2,k) = \ds\frac{9n}{10} - O(k)$, and we have the following conjecture
for the case $r = 6$.

\begin{conj}
Let $n,k \in \N$, with $n$ sufficiently large compared to $k$. Then
$$m(n,6,2,k) = \frac{3n}{4} - O(k).$$
\end{conj}

We remark that the upper bound, $m(n,6,2,k) \le \ds\frac{3n}{4}$, follows
from the construction in Lemma~\ref{s/root} below, with $R = 4$. We suspect
that the following problem is not easy.

\begin{prob}
Determine $m(n,r,2,k)$ (up to an error term depending on $r$ and $k$) for
those $r \in \N$ such that $r + 1$ is not a power of 2.\\
\end{prob}

\section{The jump at $2s = r$}\label{2s=r}

We next turn to the range $2s \approx r$, where the function $m(n,r,s,k)$
`jumps' from $(c + o(1))n$ with $c < 1$, to $n - f(k)$. We shall prove
Theorems~\ref{2ss} and \ref{2s-1}, which describe this transition quite
precisely. We begin with an instant corollary of Theorem~\ref{21k}, which
turns out to give the exact minimum when $2s = r$.

\begin{lemma}\label{2sslower}
Let $n,s,k \in \N$, with $n \ge 13k - 15$. Then $$m(n,2s,s,k) \ge n - 2k +
2.$$
\end{lemma}

\begin{proof}
Let $n,s,k \in \N$ with $n \ge 13k - 15$, let $f$ be a $(2s)$-colouring of
$E(K_n)$, and let $S \subset [2s]$ with $|S| = s$. Define the 2-colouring
$f_S$ induced by $f$ and $S$ by $f_S(e) = 1$ if $f(e) \in S$, and $f_S(e) =
2$ otherwise. By Theorem~\ref{21k}, since $n \ge 13k - 15$, there exists a
monochromatic $k$-connected subgraph $H$ of $K_n$ (in the colouring $f_S$)
on at least $n - 2k + 2$ vertices. $H$ uses at most $s$ colours in the
colouring $f$, so $m(n,2s,s,k) \ge n - 2k + 2$.
\end{proof}

Next we prove the matching upper bound. The colouring which gives the bound
is a generalization of the 2-colouring of Bollob\'as and
Gy\'arf\'as~\cite{BG}.

\begin{lemma}\label{2ssupper}
Let $n,s,k \in \N$, with $n \ge 2\ds{{2s} \choose s}(k-1) + 1$. Then
$$m(n,2s,s,k) \le n - 2k + 2.$$
\end{lemma}

\begin{proof}
Let $n,s,k \in \N$, with $n \ge 2\ds{{2s} \choose s}(k-1) + 1$. For each
subset $T \subset [2s]$ with $|T| = s$, let $A_T$ and $B_T$ be subsets of $V
= V(K_n)$ of size $k - 1$, with the sets $\{A_T, B_T : T \subset [2s], |T| =
s\}$ pairwise disjoint. Let $W = V \setminus \bigcup_T A_T \cup B_T$, so
$|W| = n - 2\ds{{2s} \choose s}(k - 1) \ge 1$. Define a $(2s)$-colouring $f$
of $E(K_n)$ as follows. Let
\begin{itemize}
\item $f(\{i,j\}) \in [2s] \setminus T$ if $i \in W$ and $j \in A_T \cup
B_T$,
\item $f(\{i,j\}) \in [2s] \setminus (T \cup T')$ if $i \in A_T \cup B_T$,
$j \in A_{T'} \cup B_{T'}$ and $T' \neq T^c$,
\end{itemize}
and for each $s$-set $T$ with $1 \notin T \subset [2s]$, let
\begin{itemize}
\item $f(\{i,j\}) = 1$ if $i \in A_T$ and $j \in A_{T^c}$, or $i \in B_T$
and $j \in B_{T^c}$, and
\item $f(\{i,j\}) \in T$ if $i \in A_T$ and $j \in B_{T^c}$, or $i \in B_T$
and $j \in A_{T^c}$.
\end{itemize}

Now, suppose $H$ is a $k$-connected subgraph of $K_n$ using at most $s$
colours; let $T \subset [2s]$, with $|T| = s$, be a fixed $s$-set containing
every colour used in $H$. We claim that $H$ contains no vertex of $A_T \cup
B_T$. Indeed, suppose $u \in V(H) \cap A_T$ say (the proof if $u \in V(H)
\cap B_T$ is identical), and let $v \in W$ (recall that $|W| \ge 1$).

Observe that since $H$ used only colours from $T$, $\Gamma_H(u) \subset
A_{T^c} \cup B_{T^c}$. Moreover, if $1 \in T$ then $\Gamma_H(u) \subset
A_{T^c}$, and if $1 \notin T$ then $\Gamma_H(u) \subset B_{T^c}$. So, if we
set $H' = H - A_{T^c}$ if $1 \in T$ and $H' = H - B_{T^c}$ if $1 \notin T$,
it is clear that $u$ and $v$ are disconnected in $H'$. Since $|A_{T^c}| =
|B_{T^c}| = k - 1$, this contradicts the assumption that $H$ is
$k$-connected. This proves the claim.

We have shown that $H$ contains no vertex of $A_T \cup B_T$. Since $|A_T
\cup B_T| = 2k - 2$, and $H$ was an arbitrary $k$-connected subgraph using
at most $s$ colours, this completes the proof of the lemma.
\end{proof}

\begin{rmk}
Note that when $s = 1$ the bound $n \ge 2\ds{{2s} \choose s}(k-1) + 1$
reduces to $n \ge 4k - 3$, and the construction reduces to that of
Bollob\'as and Gy\'arf\'as~\cite{BG}. Notice also that the construction may
be altered slightly to give the bound $m(n,2s,s,k) \le n - 2a$ if $n \ge
2\ds{{2s} \choose s}a + 1$, for each $a \le k - 1$.
\end{rmk}

By Lemma~\ref{2sslower}, any $r$-colouring of $E(K_n)$ contains a
$k$-connected subgraph, using at most $s$ colours, on at least $n - 2k + 2$
vertices, if $r \le 2s$. Suppose $s$ is decreased a little, are similar
statements are still true? In particular, for which $s$ can we always find a
$k$-connected subgraph on at least $n - g(k)$ vertices (for some function
$g$)? Or on at least $n - o(n)$ vertices? It turns out that the answer is
the same in each case: if and only if $2s \ge r$.

\begin{lemma}\label{s<2r}
Let $n,r,s,k \in \N$. If $2s < r$, then $$m(n,r,s,k) \le \left\lceil \left(1
- \ds{r \choose s}^{-1} \right)n \right\rceil.$$
\end{lemma}

\begin{proof}
Let $n,r,s,k \in \N$, with $s < 2r$. Partition $V = V(K_n)$ into $t = \ds{r
\choose s} \ge 3$ subsets $\{A_T : T \subset [r], |T| = s\}$ of near equal
size. Colour an edge between $A_T$ and $A_{T'}$ with any colour from the set
$[r] \setminus (T \cup T')$. Since $|T \cup T'| \le 2s < r$, such a colour
always exists. Colour the edges within the sets $A_T$ arbitrarily.

Let $T$ be any subset of $[r]$ with $|T| = s$. Note that $$\left \lfloor
\ds\frac{n}{t} \right \rfloor \; \le \; |A_T| \; \le \; \left \lceil
\ds\frac{n}{t} \right \rceil \; \le \; \ds\frac{n}{2},$$ since the parts are
of near equal size and $t \ge 3$. There are no edges of colour $T$ between
$A_T$ and $V \setminus A_T$, so the largest $k$-connected subgraph using the
colours of $T$ has order at most $$\max\{|A_T|,n - |A_T|\} \; \le \; n -
\left\lfloor \ds\frac{n}{t} \right\rfloor  \; = \; \left\lceil \left(1 -
\ds\frac{1}{t} \right) n \right\rceil.$$ Since $T$ was an arbitrary subset
of $[r]$ of size $s$, the lemma follows.
\end{proof}

Theorem~\ref{2ss} now follows immediately from Lemma~\ref{2sslower},
\ref{2ssupper} and \ref{s<2r}.

\begin{proof}[Proof of Theorem~\ref{2ss}]
Let $n,s,k \in \N$ with $n \ge 2\ds{{2s} \choose s}(k-1) + 1$, and $n \ge
13k - 15$. Since $n \ge 13k - 15$, we have $m(n,2s,s,k) \ge n - 2k + 2$ by
Lemma~\ref{2sslower}, and since $n \ge 2\ds{{2s} \choose s}(k-1) + 1$, we
have  $m(n,2s,s,k) \le n - 2k + 2$ by Lemma~\ref{2ssupper}. Hence
$m(n,2s,s,k) = n - 2k + 2$.

The moreover part of the theorem follows by Lemma~\ref{s<2r}.
\end{proof}

When $s = 1$ it is easy to modify the construction of Bollob\'as and
Gy\'arf\'as~\cite{BG} to give $m(n,2,1,k) = 0$ when $n \le 4k - 4$, and they
conjectured that the function jumps at this point, from $0$ to $n - 2k + 2$.
For $s \ge 2$ however, we have little clue how this transition occurs.

\begin{prob}
Determine $m(n,2s,s,k)$ when $n \le 2\ds{{2s} \choose s}(k - 1)$.
\end{prob}

Theorem~\ref{2ss} implies that $2s = r$ marks a sort of `threshold' for
$m(n,r,s,k)$: when $2s < r$ there is a constant $\eps(r,s,k) < 1$ such that
$m(n,r,s,k) < (1 - \eps(r,s,k))n$ for every (large) $n \in \N$; when $2s \ge
r$, no such constant exists, and in fact $m(n,r,s,k) \ge n - 2k + 2$ (if $n
\ge 13k$). Putting it concisely, the function `jumps' from $n - \Omega(n)$
to $n - 2k + 2$. It is natural to ask what how large $\eps(r,s,k)$ can be,
given $2s < r$. Theorem~\ref{2s-1} shows that the maximum is exactly
$\ds{{2s+1} \choose s}^{-1}$.\\

Say that a set $A \subset V$ is \emph{$(k$-$)$connected by the set $X
\subset [r]$} if the graph with vertex set $A$ and edge set $f^{-1}(X)$ is
($k$-)connected. We shall now prove Theorem~\ref{2s-1}, beginning with the
case $k = 1$.

\begin{thm}
Let $n,s \in \N$, with $s \ge 2$. Then
$$m(n,2s+1,s,1) = \left\lceil \left(1 - {{2s+1} \choose s}^{-1} \right) n
\right\rceil.$$
\end{thm}

\begin{proof}
The upper bound follows from Lemma~\ref{s<2r}; we shall prove the lower
bound. Let $n,s \in \N$, with $s \ge 2$, let $r = 2s + 1$, and let $f$ be an
$r$-colouring of $E(K_n)$. Let $\S = \{S \subset [r] : |S| = s\}$, and for
each subset $S \in \S$, let $A_S$ be a set of maximum order which is
connected by $S$, and let $B_S = V \setminus A_S$. Note that $|A_S| \ge 1$
for every $S$. We are required to show that there exists a set $S \in \S$
such that $|B_S| \le \ds{{2s+1} \choose s}^{-1} n$, and we shall do so by
showing that either $A_S$ is large for some $S \in \S$, or the sets $\{B_S :
S \in \S\}$ must be pairwise disjoint. To make things easier to follow, we
break the proof into several cases.\\

\noindent Case 1: There exist $S,T \in \S$ such that $A_S = A_T$ but $S \neq
T$.\\

If $A_S = A_T = V$ then we are done, so assume that $B_S$ is non-empty.
Since $A_S$ and $A_T$ are maximal, every edge between $A_S$ and $B_S$ must
have a colour from the set $U = \ol{S} \cap \ol{T}$. Note that $|U| = 2s + 1
- |S \cup T| \le s$, since $S \neq T$. Thus $|A_W| = n$ for any $U \subset W
\in \S$, and we are done.\\

\noindent Case 2: There exist $S,T \in \S$ such that $A_S \cap A_T =
\emptyset$.\\

Since $A_S$ and $A_T$ are maximal, every edge between $A_S$ and $A_T$ must
have some colour from $U = \ol{S} \cap \ol{T}$. Note again that $|U| \le s$,
since $S \neq T$. Now, let $B = B_S \cap B_T$, let $$C = \{v \in B : \exists
w \in A_S \cup A_T\textup{ with }f(vw) \in U\},$$ and let $D = B \setminus
C$. The set $V \setminus D$ is connected by $U$, so if $|D| = 0$ then $|A_W|
= n$ for any $U \subset W \in \S$, and we are done.

So assume that $|D| \neq 0$ and let $u \in D$. Since $u \notin A_S$, edges
between $u$ and $A_S$ do not have colours from $S$. Also, by the definition
of $D$, these edges do not have colours from $U$. Thus they must have
colours from $[r] \setminus (S \cup U) \subset T$. Similarly, edges between
$u$ and $A_T$ must have colours from $S$. This is true for any vertex in
$D$, so the set $A_S \cup D$ is connected by $T$, and the set $A_T \cup D$
is connected by $S$. But $A_S$ and $A_T$ have maximum order, so $|A_S| + |D|
\le |A_T|$, and $|A_T| + |D| \le |A_S$, which implies that $|D| = 0$, a
contradiction.\\

\noindent Case 3: There exist $S,T \in \S$ such that $A_S \cap A_T \neq
\emptyset$, $A_S \not\subset A_T \not\subset A_S$, $A_S \cup A_T \neq V$ and
$|S \cup T| \ge s + 2$.\\

Let $u \in B = V \setminus (A_S \cup A_T)$, and let $C = A_S \cap A_T \neq
\emptyset$. Since $u \notin A_S$ and $u \notin A_T$, the edges between $u$
and $C$ must all have colours from $U = \ol{S} \cap \ol{T}$. Similarly, all
edges between $A_S \setminus A_T$ and $A_T \setminus A_S$ must have colours
from $U$. Thus, the entire vertex set $V$ is connected by $U \cup \{i\}$,
where $f(vw) = i$ for some $v \in B \cup C$ and $w \in A_S \triangle A_T$.
Now, since $|S \cup T| \ge s + 2$, it follows that $|U| \le s - 1$, so $|U
\cup \{i\}| \le s$ and we are done.\\

\noindent Case 4: There exist $S,T \in \S$ such that $A_S \cap A_T \neq
\emptyset$, $A_S \not\subset A_T \not\subset A_S$, $A_S \cup A_T \neq V$ and
$|S \cup T| = s + 1$.\\

We shall show that Case 3 still holds. Indeed, let $B$, $C$ and $U$ be as in
Case 3, and note that $|U| = s$, and that $U \cap S = U \cap T = \emptyset$.
As before, the sets $B \cup C$ and $A_S \triangle A_T$ are each connected by
$U$, and these sets partition the vertex set $V$. Hence, either $A_U = V$,
in which case we are done, or $A_U = B \cup C$, or $A_U = A_S \triangle
A_T$. It is simple to check that in either of the latter two cases $S$ and
$U$ satisfy the conditions of Case 3, and so we are done as before. Note in
particular that $|S \cup U| = 2s \ge s + 2$ since $s \ge 2$.\\

\noindent Case 5: There exist $S,T \in \S$ such that $A_S \supset A_T$ but
$S \neq T$.\\

This case is a little more complicated than the first four. First we shall
show that $A_R \supset A_T$ for \emph{every} set $R \in \S$.

If $A_S = V$ we are done, so we may assume that $B_S$ is non-empty. Since
$A_S$ and $A_T$ are maximal, all edges between $B_S$ and $A_T$ must have
colours from $U = \ol{S} \cap \ol{T}$. Choose $W \in \S$ such that $U
\subset W$, and let $C_W$ be the maximal set connected by $W$ containing
$B_S \cup A_T$. We shall show that $|C_W| > n/2$, and so $C_W = A_W$. If
$C_W = V$ then we are done, so assume that $D_W = V \setminus C_W$ is
non-empty. Now, no edge between $D_W \subset A_S$ and $B_S \subset C_W$ can
have a colour from $S$ or from $U$, so all of these edges have colours from
$T$. But $A_T$ was chosen to have maximal size, so $|A_T| \ge |D_W| +
|B_S|$. Hence $$|C_W| \; \ge \; |A_T| + |B_S| \; \ge \; |D_W| + 2|B_S| \; >
\; |D_W|,$$ and so $C_W = A_W$ as claimed.

We have shown that $A_W \supset A_T$ for every $\ol{S} \cap \ol{T} \subset W
\in \S$, so in particular we can choose $W$ so that $T \cap W = \emptyset$.
Let $\{i\} = [r] \setminus (T \cup W)$. Now, by the method of the previous
paragraph, $A_R \supset A_T$ for any $R \in \S$ with $i \in R$. In
particular, if $X = W \triangle \{i,j\}$ with $j \in W$, then $A_X \supset
A_T$. Once again applying the method of the previous paragraph, we infer
that $A_R \supset A_T$ for any $R \in \S$ with $j \in R$. Since $j$ was an
arbitrary member of $W = [r] \setminus (T \cup \{i\})$, we have proved that
$A_R \supset A_T$ for \emph{every} set $R \in \S$, as claimed.

We next claim that $B_Q \cap B_R = \emptyset$ for every $Q, R \in \S
\setminus \{T\}$ with $Q \neq R$. Indeed, if $B_Q \cap B_R \neq \emptyset$
and $B_Q \not\subset B_R \not\subset B_Q$, then we are in either Case 3 or
Case 4, since $B_Q \cup B_R \subset B_T \neq V$. But if, on the other hand,
$B_Q \subset B_R$ say, then $A_Q \supset A_R$, so $A_P \supset A_R$ for
every $P \in \S$ as above, and in particular $A_T \supset A_R$. But then
$A_R = A_T$, and we are in Case 1. Hence $B_Q \cap B_R = \emptyset$ for
every $Q, R \in \S \setminus \{T\}$ with $Q \neq R$, as claimed.

Now, simply observe that for any pair $W, X \in \S$ such that $\ol{T}
\subset W, X$ and $W \neq X$, every edge between $B_W$ and $B_X$ must have a
colour from $T$, since $[r] \setminus (W \cup X) \subset T$. Since $B_W$ and
$B_X$ are disjoint, and $A_T$ was chosen to be maximal, it follows that
$|A_T| \ge |B_W| + |B_X|$.

Hence, recalling that $A_T \cap B_R = B_R \cap B_Q = \emptyset$ for every
$Q,R \in \S \setminus \{T\}$ with $Q \neq R$, we obtain
$$n \; \ge \; \left( \sum_{T \neq R \in \S} |B_R| \right) + |B_W| + |B_X| \;
\ge \; \left( {{2s+1} \choose s} + 1 \right) \min_{R \in \S} |B_R|,$$ and
thus $\ds\min_{R \in \S} |B_R| \le \left( {{2s+1} \choose s} + 1
\right)^{-1}n$, as required.\\[+2ex]

Finally, suppose that none of Cases 1--5 hold. The only remaining
possibility is that $A_S \cup A_T = V$ for every pair $S,T \in \S$ with $S
\neq T$, and therefore that $B_S \cap B_T = \emptyset$ for every such pair.
But now we have
$$n \; \ge \; \sum_{R \in \S} |B_R| \; \ge \; {{2s+1} \choose s} \min_{R \in
\S} |B_R|,$$ and so $\ds\min_{R \in \S} |B_R| \le {{2s+1} \choose s}^{-1}n$,
and we are done.
\end{proof}

The proof for general $k$ is similar, but we shall need some of the tools
from Section~\ref{tools}: to be precise, we shall use Lemmas~\ref{allSbig},
\ref{21ktech} and \ref{2s-1bip}.

\begin{proof}[Proof of Theorem~\ref{2s-1}]
The upper bound again follows from Lemma~\ref{s<2r}; we shall prove the
lower bound. Let $n,s,k \in \N$, with $s \ge 2$ and $n \ge 100\ds{{2s + 1}
\choose s}^2k^2$. Let $r = 2s + 1$, and let $f$ be an $r$-colouring of
$E(K_n)$. Let $\S = \{S \subset [r] : |S| = s\}$ as before, and for each
subset $S \in \S$, let $A_S$ be a set of maximum order which is
$k$-connected by $S$, and let $B_S = V \setminus A_S$. We are required to
show that there exists a set $S \in \S$ such that $|B_S| \le \ds{{2s+1}
\choose s}^{-1} n + 2{{2s + 1} \choose s}k$.

Let us assume, for a contradiction, that no such set $S$ exists. We begin by
using Lemma~\ref{allSbig} to show that $|A_S| \ge (6s + 78)k$ for
\emph{every} set $S \in \S$. Indeed, let $S \in \S$, and let $T,U \subset
[r]$ satisfy $S \cup T \cup U = [r]$. We have $$q_k(W) \; \le \; \left( 1 -
\ds{{2s+1} \choose s}^{-1}\right) n - 2{{2s + 1} \choose s}k \; < \; n -
k\sqrt{n}$$ for each $W \in \{T,U\}$, and $n \ge 25k^2$, so by
Lemma~\ref{allSbig},
\begin{eqnarray*}
q_k(S) & \ge & n \: - \: 9k\sqrt{n} \: - \: q_k(T) \\[+1ex]
& \ge & \ds{{2s+1} \choose s}^{-1}n \: + \: 2{{2s + 1} \choose s}k \: - \:
9k\sqrt{n}\\[+1ex]
& = & \ds{{2s+1} \choose s}^{-1} \sqrt{n} \left( \sqrt{n} - 9k\ds{{2s+1}
\choose s} \right) \: + \: 2{{2s + 1} \choose s}k \\[+1ex]
& \ge & 10{{2s + 1} \choose s}k^2 \: + \: 2{{2s + 1} \choose s}k \; > \; (6s
+ 78)k
\end{eqnarray*}
since $\sqrt{n} \ge 10k\ds{{2s + 1} \choose s}$. Thus $|A_S| \ge (6s + 78)k$
for every $S \in \S$, as claimed. We must once again consider five cases.\\

\noindent Case 1: There exist $S,T \in \S$ such that $|A_S \triangle A_T|
\le (4s + 62)k$.\\

If $|V \setminus (A_S \cup A_T)| \le 2k$, then $|A_S| \ge n - (4s + 64)k$
and we are done, so assume that $|B| \ge 2k$, where $B = B_S \cap B_T$.
Since $A_S$ is maximal, a vertex $v \in B$ can send at most $k-1$ edges with
colours from $S$ into $A_S$, and similarly $v$ can send at most $k - 1$
edges with colours from $T$ into $A_T$. Therefore, each vertex of $B$ sends
at most $2k - 2$ edges with colours not in the set $U = \ol{S} \cap \ol{T}$
into $A = A_S \cap A_T$.

Now, simply apply Lemma~\ref{21ktech} to the bipartite graph $G$ with parts
$A$ and $B$, and edges with colours from the set $U$, and with $a = b = 2k$.
Note that $$|A| \; = \; \frac{|A_S| + |A_T| - |A_S \triangle A_T|}{2} \; \ge
\; 9k \; = \; 4b + k,$$ and $|B| \ge 2k$, so $G$ contains a $k$-connected
subgraph on more than $|G| - 4k$ vertices. But this subgraph uses only
colours from $U$, and $|G| = n - |A_S \triangle A_T| \ge n - (4s + 62)k$, so
in this case $q_k(U) \ge n - (4s + 66)k$ and we are done.\\

\noindent Case 2: There exist $S,T \in \S$ such that $|A_S \cap A_T| \le (2s
+ 16)k$.\\

Since $A_S$ and $A_T$ are maximal, each vertex of $A_S \setminus A_T$ sends
at most $k - 1$ edges with colours from $T$ into $A_T \setminus A_S$, and
similarly each vertex of $A_T \setminus A_S$ sends at most $k - 1$ edges
with colours from $S$ into $A_S \setminus A_T$. We also have $|A_S \setminus
A_T| = |A_S| - |A_S \cap A_T| \ge 15k$, and similarly $|A_T \setminus A_S|
\ge 15k$, so we may apply Lemma~\ref{2s-1bip} to obtain a set $X \subset A_S
\triangle A_T$ with $|X| \ge |A_S \triangle A_T| - 7k$, which is
$k$-connected by $U = \ol{S} \cap \ol{T}$.

Let $B = B_S \cap B_T$, let $$C = \{v \in B : |\{w \in X : f(vw) \in U\}|
\ge k\},$$ and let $D = B \setminus C$. Note that the set $X \cup C$ is
$k$-connected by $U$, so $|A_W| \ge n - 12k - |D|$ for any $U \subset W \in
\S$. Therefore we may assume that $|D| > (4s + 50)k$.

Now, each vertex $u \in D$ sends at most $k - 1$ edges with colours from $S$
into $A_S$ (since $u \notin A_S$), and at most $k - 1$ edges with colours
from $U$ into $X$ (by the definition of $D$). Apply Lemma~\ref{21ktech} to
the bipartite graph with parts $D$ and $X \cap A_S$, and edges with colours
from $T$, and with $a = b = 2k$. Note that $|D| \ge 2k$ and $$|X \cap A_S|
\; \ge \; |A_S| - |A_S \cap A_T| - 5k \; \ge \; 9k \; = \; 4b + k,$$ so the
conditions of the lemma hold. Thus there exists a $k$-connected subgraph of
$K_n$, using colours only from $T$, with at least $$|D| + |X \cap A_S| - 4k
\; \ge \; |D| + |A_S| - (2s + 25)k$$ vertices. Similarly, there exists a
$k$-connected subgraph of $K_n$, using colours only from $S$, with at least
$|D| + |A_T| - (2s + 25)k$ vertices. But $A_S$ and $A_T$ have maximum order,
so $|A_S| + |D| - (2s + 25)k \le |A_T|$, and $|A_T| + |D| - (2s + 25)k \le
|A_S|$, which implies that $|D| \le (4s + 50)k$, a contradiction.\\

\noindent Case 3: There exist $S,T \in \S$ such that $|A_S \cap A_T| \ge (2s
+ 5)k$, $|A_S \setminus A_T|, |A_T \setminus A_S| \ge (s + 13)k$, $|V
\setminus (A_S \cup A_T)| \ge 2k$ and $|S \cup T| \ge s + 2$.\\

Let $B = V \setminus (A_S \cup A_T)$, $C = A_S \cap A_T$, $X = A_S \setminus
A_T$ and $Y = A_T \setminus A_S$, so $|B| \ge 2k, |C| \ge (2s + 5)k \ge 9k$,
and $|X|, |Y| \ge (s + 13)k \ge 15k$. Since $A_S$ and $A_T$ are maximal, a
vertex of $B$ can send at most $2k - 2$ edges with colours from $S \cup T$
into $C$, thus by Lemma~\ref{21ktech} there exists a $k$-connected subgraph
$H_1$ of $K_n[B \cup C]$, using only colours of $U = \ol{S} \cap \ol{T}$, on
at least $|B| + |C| - 4k$ vertices. Furthermore, each vertex of $X$ sends at
most $k-1$ edges with colours from $T$ into $Y$, and each vertex of $Y$
sends at most $k-1$ edges with colours from $S$ into $X$, so by
Lemma~\ref{2s-1bip} there exists a $k$-connected subgraph $H_2$ of $K_n[X
\cup Y]$, using only colours of $U = \ol{S} \cap \ol{T}$, on at least $|X| +
|Y| - 7k$ vertices. Thus $|V(H_1) \cup V(H_2)| \ge n - 11k$.

It remains only to `$k$-connect' $H_1$ and $H_2$ using
Observation~\ref{linkup}. To be precise, for each $i \in S \cup T$ let
$$D_i = \{v \in V(H_1) : |\{w \in V(H_2) : f(vw) \in U \cup \{i\}\}| \ge
k\},$$
and note that $\bigcup_i D_i = V(H_1)$, since $|H_2| \ge |X| + |Y| - 7k \ge
2sk$. Choose $j \in S \cup T$ such that $|D_j| \ge |H_1|/2s$, and note that
$|H_1| \ge |B| + |C| - 4k \ge 2sk$, so $|D_j| \ge k$.

Now, by Observation~\ref{linkup}, $V(H_1) \cup V(H_2)$ is $k$-connected by
$U \cup \{j\}$, and so $q_k(U \cup \{j\}) \ge n - 11k$. But $|S \cup T| \ge
s + 2$, so $|U \cup \{j\}| \le s$, and we are done.\\

\noindent Case 4: There exist $S,T \in \S$ such that $|A_S \cap A_T| \ge (2s
+ 16)k$, $|A_S \setminus A_T|, |A_T \setminus A_S| \ge (2s + 19)k$, $|V
\setminus (A_S \cup A_T)| \ge (s + 17)k$ and $|S \cup T| = s + 1$.\\

We shall show that either Case 3 still holds, or $|A_U| \ge n - 11k$, where
as usual $U = \ol{S} \cap \ol{T}$ (note that now $|U| = s$, since $|S \cup
T| = s + 1$). Let $B$, $C$, $X$ and $Y$ be as in Case 3, so $|B| \ge (s +
17)k$, $|C| \ge (2s + 16)k$, and $|X|, |Y| \ge (2s + 19)k$. As before, we
can find disjoint $k$-connected subgraphs $H_1$ and $H_2$, which use only
edges with colours from $U$, such that $V(H_1) \subset B \cup C$, $V(H_2)
\subset X \cup Y$, $|V(H_1)| \ge |B| + |C| - 4k$ and $|V(H_2)| \ge |X| + |Y|
- 7k$.

Consider $A_U$. Clearly $|A_U| \ge \max\{|H_1|, |H_2|\} \ge (n - 11k)/2 \ge
13k$, but any set of order $13k$ intersects either $V(H_1)$ or $V(H_2)$ (or
both) in at least $k$ vertices. Hence, by Observation~\ref{GcupH}, either
$A_U \supset V(H_1)$ and $|A_U \cap V(H_2)| \le k - 1$, or $A_U \supset
V(H_2)$ and $|A_U \cap V(H_2)| \le k - 1$, or $A_U \supset V(H_1) \cup
V(H_2)$. In the third case we have $|A_U| \ge n - 11k$; we claim that in
either of the first two (sub)cases, $S$ and $U$ satisfy the conditions of
Case 3.\\

\noindent Subcase (a): If $A_U \supset V(H_1)$ and $|A_U \cap V(H_2)| \le k
- 1$, then
\begin{eqnarray*}
|A_S \cap A_U| & \ge & |V(H_1) \cap C| \; \ge \; |C| - 4k \; \ge \; (2s +
5)k,\\
|A_S \setminus A_U| & \ge & |V(H_2) \cap X| - k \; \ge \; |X| - 8k \; \ge \;
(s+13)k,\\
|A_U \setminus A_S| & \ge & |V(H_1) \cap B| \; \ge \; |B| - 4k \; \ge \; (s
+ 13)k,\textup{ and }\\
|V \setminus (A_S \cup A_U)| & \ge & |V(H_2) \cap Y| - k \; \ge \; |Y| - 8k
\; \ge \; 2k.\\[-1ex]
\end{eqnarray*}

\noindent Subcase (b): Similarly, if $A_U \supset V(H_2)$ and $|A_U \cap
V(H_1)| \le k - 1$, then
\begin{eqnarray*}
|A_S \cap A_U| & \ge & |V(H_2) \cap X| \; \ge \; |X| - 7k \; \ge \; (2s +
5)k,\\
|A_S \setminus A_U| & \ge & |V(H_1) \cap C| - k \; \ge \; |C| - 5k \; \ge \;
(s+13)k,\\
|A_U \setminus A_S| & \ge & |V(H_2) \cap Y| \; \ge \; |Y| - 7k \; \ge \; (s
+ 13)k,\textup{ and}\\
|V \setminus (A_S \cup A_U)| & \ge & |V(H_1) \cap B| - k \; \ge \; |B| - 5k
\; \ge \; 2k. \\[-1ex]
\end{eqnarray*}
Also $|S \cup U| = 2s \ge s + 2$ since $s \ge 2$, and so $S$ and $U$ satisfy
the conditions of Case 3, as claimed. Hence we are done as in that case.\\

\noindent Case 5: There exist $S,T \in \S$ such that $|A_T \setminus A_S|
\le (2s + 19)k$ but $S \neq T$.\\

This case is once again a little more complicated than the first four. First
we shall show that $|A_T \setminus A_R| \le (2s + 31)k$ for \emph{every} set
$R \in \S$.

Let $B = B_S \cap B_T$ and $C = A_S \cap A_T$, as in Cases 3 and 4. If
$|A_S| \ge n - (4s + 48)k$ we are done, so we may assume that $|B_S| \ge (4s
+ 48)k$, and hence that
$$|B| \; = \; |B_S \cap B_T| \; = \; |B_S| - |A_T \setminus A_S| \; \ge \;
(2s + 29)k.$$
Also recall that $|A_T| \ge (6s + 78)k$, so
$$|C| \; = \; |A_S \cap A_T| \; = \; |A_T| - |A_T \setminus A_S| \; \ge \;
9k.$$

Now, since $A_S$ and $A_T$ are maximal, each vertex of $B$ sends at most $2k
- 2$ edges with colours from $S \cup T$ into $C$. Let $U = \ol{S} \cap
\ol{T}$, and apply Lemma~\ref{21ktech} in the usual way (with $a = b = 2k$)
to obtain a $k$-connected subgraph $H$ of $K_n[B \cup C]$, using only
colours from the set $U$, on at least $|B| + |C| - 4k$ vertices. Choose $W
\in \S$ such that $U \subset W$, and let $C_W$ be a set of maximum size,
containing $V(H)$, which is $k$-connected by $W$. We shall show that $|C_W|
\ge (n + k)/2$, and deduce that $C_W = A_W$.

Indeed, we are done if $|C_W| \ge n - (2s + 38)k$, since $|A_W| \ge |C_W|$,
so let $D_W = V \setminus C_W$ and assume that $|D_W| \ge (2s + 38)k$. Then
\begin{eqnarray*}
|D_W \cap A_S| & \ge & |D_W| - |A_T \setminus A_S| - |B \cap D_W|\\[+1ex]
& \ge & |D_W| - (2s + 19)k - 4k \; \ge \; 15k,
\end{eqnarray*}
and also
$$|B_S \cap V(H)| \; \ge \; |B| - 4k \; \ge \; 15k.$$
Now, a vertex of $D_W$ can send at most $k - 1$ edges with colours from $W$
into $V(H)$, and a vertex of $B_S$ can send at most $k - 1$ edges with
colours from $S$ into $A_S$. Therefore, by Lemma~\ref{2s-1bip}, there exists
a $k$-connected subgraph of $(D_W \cap A_S) \cup (B_S \cap V(H))$, using
only colours from the set $T$, on at least
$$ |D_W \cap A_S| + |B_S \cap V(H)| - 7k \; \ge \; |D_W| + |B| - (2s +
34)k$$
vertices. But $A_T$ was chosen to have maximal size, so
$$|A_T| \; \ge \; |D_W| + |B| - (2s + 34)k.$$
Hence
\begin{eqnarray*}
|C_W| & \ge & |B| + |C| - 4k \; = \; |B| + |A_T| - |A_T \setminus A_S| -
4k\\[+1ex]
& \ge & |D_W| + 2|B| - (4s + 57)k \; \ge \; |D_W| + k,
\end{eqnarray*}
since $|B| \ge (2s + 29)k$. Therefore $|C_W| \ge (n + k)/2$, as claimed. But
now any subset of $V$ of size at least $|C_W|$ must intersect $C_W$ in at
least $k$ vertices, and so by Observation~\ref{GcupH}, any $k$-connected
subgraph on at least $|C_W|$ vertices must contain $C_W$. In particular,
$A_W \supset C_W$, since $|A_W| \ge |C_W|$ by definition, But $C_W$ was
chosen to have maximum size, so we must have $C_W = A_W$.

Since $V(H) \subset A_W$, we have shown that
$$|A_T \setminus A_W| \; \le \; |A_T \setminus A_S| + |C \setminus V(H)| \;
\le \; (2s + 23)k$$ for every set $W \in \S$ with $\ol{S} \cap \ol{T}
\subset W \in \S$. In particular, we may choose $W$ so that $T \cap W =
\emptyset$. Now applying the method of the previous paragraphs to the sets
$T$ and $W$, we deduce that $|A_T \setminus A_R| \le (2s + 27)k$ for any $R
\in \S$ with $i \in R$, where $\{i\} = \ol{T} \cap \ol{W}$. In particular,
if $X = W \triangle \{i,j\}$ with $j \in W$, then $|A_T \setminus A_X| \le
(2s + 27)k$. Once again applying the method of the previous paragraphs, we
infer that $|A_T \setminus A_R| \le (2s + 31)k$ for any $R \in \S$ with $j
\in R$. Since $j$ was an arbitrary member of $[r] \setminus (T \cup \{i\})$,
we have proved that $|A_T \setminus A_R| \le (2s + 31)k$ for \emph{every}
set $R \in \S$, as claimed.

We shall next show that either we are in Case 1, 3 or 4, or $|B_Q \cap B_R|
< (s + 17)k$ for every $Q, R \in \S \setminus \{T\}$ such that $Q \neq R$,
and moreover $|B_Q \cap B_R| < 2k$ if $|Q \cup R| \ge s + 2$. Indeed, let
$Q, R \in \S \setminus \{T\}$ with $Q \neq R$, and let $|B_Q \cap B_R| \ge
2k$. Suppose first that  $|B_Q \setminus B_R| \le (2s + 19)k$. Then $|A_R
\setminus A_Q| \le (2s + 19)k$, and so $|A_R \setminus A_P| \le (2s + 31)k$
for every $P \in \S$, as above, and in particular $|A_R \setminus A_T| \le
(2s + 31)k$. But now $|A_R \triangle A_T| \le (4s + 62)k$, and we are in
Case 1.

So suppose next that $|B_Q \setminus B_R|, |B_R \setminus B_Q| \ge (2s +
19)k$. Note that
\begin{eqnarray*}
|A_Q \cap A_R| & \ge & |A_T| - |A_T \setminus A_R| - |A_T \setminus
A_Q|\\[+1ex]
& \ge & |A_T| - 2(2s + 31)k \; \ge \; (2s + 16)k,
\end{eqnarray*}
since $|A_T| \ge (6s + 78)k$, that $|A_Q \setminus A_R|, |A_R \setminus A_Q|
\ge (2s + 19)k$, and that
$$|V \setminus (A_Q \cup A_R)| \; = \; |B_Q \cap B_R|.$$
Thus if $|Q \cup R| \ge s + 2$ we are in Case 3, and if $|Q \cup R| = s + 1$
and $|B_Q \cap B_R| \ge (s + 17)k$ then we are in Case 4. Hence either we
are done as before, or $|B_Q \cap B_R| < (s + 17)k$ for every $Q, R \in \S
\setminus \{T\}$ with $Q \neq R$, and moreover $|B_Q \cap B_R| < 2k$ if $|Q
\cup R| \ge s + 2$, as claimed.

Now, let $W$ and $X$ be as described above, so $[r] \setminus (W \cup X)
\subset T$, and observe that a vertex of $B_W \setminus B_X$ sends at most
$k-1$ edges with colours from $W$ into $B_X \setminus B_W$, and similarly a
vertex of $B_X \setminus B_W$ sends at most $k-1$ edges with colours from
$X$ into $B_W \setminus B_X$. Note also that $|B_W|, |B_X| \ge (s + 32)k$,
else we are done, so
$$|B_W \setminus B_X| \; \ge \; |B_W| - |B_W \cap B_X| \; \ge \; 15k,$$ and
similarly $|B_X \setminus B_W| \ge 15k$. Hence we may apply
Lemma~\ref{2s-1bip} to the bipartite graph with parts $B_W \setminus B_X$
and $B_X \setminus B_W$ to obtain a $k$-connected subgraph on at least
$$|B_W \triangle B_X| - 7k \; \ge \; |B_W| + |B_X| - (s + 24)k$$ vertices,
using only colours from the set $T$.

Since $A_T$ was chosen to be maximal, we have $$|A_T| \; \ge \; |B_W| +
|B_X| - (s + 24)k.$$ Now, recall that for every $Q,R \in \S \setminus
\{T\}$, we have $$|B_R \setminus A_T| \; = \; |B_R| - |A_T \setminus A_R| \;
\ge \; |B_R| - (2s + 31)k,$$ $|B_Q \cap B_R| \le 2k$ if $|Q \cup R| \ge s +
2$ and $|B_Q \cap B_R| \le (s + 17)k$ if $|Q \cup R| = s + 1$. There are
exactly $\ds{{s + 1} \choose 2}$ pairs $Q,R \in \S$ such that $|Q \cap R| =
s + 1$. Thus, by inclusion-exclusion, we obtain
\begin{eqnarray*}
n & \ge & |A_T| \; + \;  \sum_{T \neq R \in \S} |B_R \setminus A_T| \; -
\sum_{\substack{Q,R \in \S \setminus \{T\},\\ Q \neq R}} |B_Q \cap
B_R|\\[+1ex]
& > & |B_W| \; + \; |B_X| \; - \; (s + 24)k \; + \; \sum_{T \neq R \in \S}
\Big( |B_R| - (2s + 31)k \Big)\\[+1ex]
& & \hspace{4cm}  - \; {{2s + 1} \choose s}^2 k \; - \; {{s + 1} \choose
2}(s + 17)k\\[+1ex]
& > & \left( {{2s+1} \choose s} + 1 \right) \min_{R \in \S} |B_R| \; - \;
\left( {{2s + 1} \choose s}^2 + {{2s + 1} \choose s}(2s + 40) \right)k.
\end{eqnarray*} Now,
\begin{eqnarray*}
\frac{n + \left( \ds{{2s + 1} \choose s}^2 + {{2s + 1} \choose s}(2s + 40)
\right)k}{\ds{{2s+1} \choose s} + 1}
& \le & \frac{n \: + \: 2\ds{{2s + 1} \choose s}^2k}{\ds{{2s+1} \choose s}}
\end{eqnarray*}
reduces to $n \ge  \ds{{2s + 1} \choose s}^2 \left(2s + 38 - \ds{{2s + 1}
\choose s} \right)k$, which is true, so
$$\ds\min_{R \in \S} |B_R| \; < \; \ds{{2s+1} \choose s}^{-1}n \: + \: 2{{2s
+ 1} \choose s}k,$$
as required.\\

Finally, suppose that none of Cases 1--5 hold. The only remaining
possibility is that $|A_S \cup A_T| \ge n - 2k$ for every pair $S,T \in \S$
with $|S \cup T| \ge s + 2$, and $|A_S \cup A_T| \ge n - (s + 17)k$ for
every pair $S,T \in \S$ with $|S \cup T| = s + 1$. But $|B_S \cap B_T| = n -
|A_S \cup A_T|$, so we have
\begin{eqnarray*}
n & \ge & \sum_{R \in \S} |B_R| \; - \sum_{\substack{Q,R \in \S \setminus
\{T\},\\ Q \neq R}} |B_Q \cap B_R|\\
& \ge & {{2s+1} \choose s} \min_{R \in \S} |B_R| \; - \; {{2s + 1} \choose
s}^2 k \; - \; {{s + 1} \choose 2}(s + 17)k,
\end{eqnarray*}
So $\ds\min_{R \in \S} |B_R| \le {{2s+1} \choose s}^{-1}n + 2{{2s + 1}
\choose s}k$, and we are done.
\end{proof}

Setting $s = 2$ we obtain the following corollary.

\begin{cor}
Let $n,k \in \N$ with $n \ge (100k)^2$. Then
$$\frac{9n}{10} \: - \: 20k \; \le \; m(n,5,2,k) \; \le \; \frac{9n}{10} \:
+ \: 1.$$
\end{cor}

\begin{rmk}
In fact one can do a little better in both directions. By taking a little
more care in the proof of Theorem~\ref{2s-1}, one easily obtains
$$m(n,5,2,k) \ge \frac{9n - 157k}{10},$$ and a simple modification of the
construction in Lemma~\ref{s<2r} gives $$m(n,5,2,k) \le \frac{9n - k +
1}{10},$$ and more generally $$m(n,2s + 1,s,k) \le \left(1 - {{2s+1} \choose
s}^{-1}\right)n - {{2s+1} \choose s}^{-1}(k - 1).$$ We are sure that neither
of these bounds is sharp.\\
\end{rmk}

\section{The jump at $s = \Theta(\sqrt{r})$}\label{s=rootr}

Perhaps the most basic question one can ask about the function $m(n,r,s,k)$
is the following: for which values of $s$ is $m(n,r,s,k)$ close to $0$, and
for which is it close to $1$? Theorem~\ref{jump} gives an asymptotic answer
to this question. We begin with an easy lemma, which gives the upper bound
in the theorem.

\begin{lemma}\label{s/root}
For every $n,r,s,k \in \N$, we have $$m(n,r,s,k) \; \le \; (s+1) \left\lceil
\ds\frac{n}{ \left\lfloor \sqrt{2r} \right\rfloor} \right\rceil.$$
\end{lemma}

\begin{proof}
Let $n,r,s,k \in \N$, let $V = V(K_n)$, and partition $V$ into $R = \lfloor
\sqrt{2r} \rfloor$ sets $V_1, \ldots, V_R$, each of size either $N = \lceil
n/R \rceil$ or $N - 1$. Noting that $\ds{R \choose 2} < r$, assign to each
pair $\{i,j\} \subset [R]$ a distinct colour $c(\{i,j\}) \in [r]$.

Let $f$ be the following $r$-colouring of $E(K_n)$: if $x \in V_i$ and $y
\in V_j$, and $i \neq j$, then set $f(xy) = c(\{i,j\})$. If $x,y \in V_i$,
then $f(xy)$ may be chosen arbitrarily. Thus $f$ is a `blow-up' of a
completely multicoloured complete graph.

Now, let $S \subset [r]$ be any subset of size $s$, and let $G$ be the
subgraph of $K_n$ with vertex set $V$ and edge set $f^{-1}(S)$. Each
component of $G$ intersects at most $s+1$ of the sets $\{V_j : j \in [R]\}$,
so since $S$ was chosen arbitrarily, we have $m(n,r,s,k) \le M(f,n,r,s,1)
\le (s+1)N$.
\end{proof}

Lemma~\ref{s/root} shows that if $s \ll \sqrt{r}$, then
$\ds\frac{m(n,r,s,k)}{n} \to 0$ as $r \to \infty$. Somewhat surprisingly,
this simple construction turns out to be asymptotically optimal. Once again,
we begin with the case $k = 1$, and prove a slightly stronger result.

\begin{thm}\label{rootk=1}
Let $n,r,s \in \N$. Then $$m(n,r,s,1) \ge \left(1 - e^{-s^2/3r} \right)n.$$
\end{thm}

\begin{proof}
Let $n,r,s \in \N$, and let $f$ be an $r$-colouring of the edges of $K_n$.
If $s = 1$ then the result is trivial, since $e^{-x} > 1 - x$ if $x > 0$,
and $m(n,r,1,1) \ge n/r$ (consider the largest monochromatic star centred at
any vertex). So let $s \ge 2$, and assume the result holds for all smaller
values of $s$. Let $t \in [s - 1]$ (we shall eventually set $t = \lceil s/2
\rceil$, but we shall delay making this choice until it is clear why it is
optimal), and let $G$ be a connected subgraph of $K_n$, using at most $t$
colours, of maximum order. Let $V = V(K_n)$, $A = V(G)$, $B = V \setminus
A$, and $T = f(E(G))$, the set of colours used by $G$. Thus (assuming $|A| <
n$), $|T| = t$. By the induction hypothesis, $|A| \ge \left(1 - e^{-t^2/3r}
\right)n$.

Now, each vertex in $B$ must send at least $t + 1$ different colours into
$A$, as otherwise the star centred at that vertex would be a connected
component, using at most $t$ colours, larger than $G$. Also, a vertex of $B$
sends no edges with colours from $T$ into $A$, since $G$ was chosen to be
maximal. For each vertex $v \in B$, choose a list $L(v)$ of $t + 1$ colours
$\{\ell_1, \ldots, \ell_{t+1}\} \subset [r] \setminus T$ which it sends into
$A$. So for each $v \in B$ and $\ell \in L(v)$, there exists a vertex $u \in
A$ such that $f(uv) = \ell$.

Let $\eps > 0$, and let $\T = \{S \subset [r] \setminus T : |S| = s - t\}$.
Suppose that $m(n,r,s,1) \le n - \eps |B|$. This means that for every set $S
\in \T$, the largest connected component in $K_n$, using only the colours $S
\cup T$, and containing $G$, avoids at least $\eps |B|$ vertices of $B$.
Hence, for each $S \in \T$ there are at least $\eps |B|$ vertices $v \in B$
such that $S \cap L(v) = \emptyset$. For each $S \in \T$, let $M(S) = \{v
\in B : S \cap L(v) = \emptyset\}$.

Now, observe that for each vertex $v \in B$, there are exactly $\ds{{r - 2t
- 1} \choose {s - t}}$ sets $S \in \T$ with $v \in M(S)$. So, summing over
$\T$, we obtain
$$\eps |B| \ds{{r-t} \choose {s-t}} \; \le \; \sum_{S \in \T} |M(S)| \; = \;
\sum_{v \in B} \sum_{S \in \T} I[v \in M(S)] \; = \; |B| \ds{{r-2t-1}
\choose {s-t}},$$ where $I[T]$ denotes the indicator function of the event
$T$, and therefore
\begin{eqnarray*}
\eps & \le & \frac{(r - 2t - 1)! \: (r - s)!}{(r - t)! \: (r - s - t - 1)!}
\; = \; \frac{(r - 2t - 1)}{(r - t)}  \: \cdots \: \frac{(r - s - t)}{(r - s
+ 1)}\\
& \le & \left( \ds\frac{r - 2t - 1}{r - t} \right)^{s-t} \; < \;
\operatorname{exp} \left( \ds\frac{-(t + 1)(s - t)}{r - t} \right).
\end{eqnarray*}
Now, set $t = \lceil s/2 \rceil$ to (approximately) maximize $\ds\frac{(t +
1)(s - t)}{r - t}$, and note that $\ds\frac{(\lceil s/2 \rceil + 1) \lfloor
s/2 \rfloor}{r - \lceil s/2 \rceil} > \ds\frac{s^2}{4r}$. Recalling that
$|B| \le e^{-t^2/3r}n \le e^{-s^2/12r}n$, we obtain
$$\eps|B| \; \le \; e^{-s^2/4r - s^2/12r}n \; = \; e^{-s^2/3r},$$ since $s
\ge 2$. Hence
$$M(f,n,r,s,1) \; \ge \; |A| + (1 - \eps)|B| \; = \; n - \eps|B| \; \ge \;
\left(1 - e^{-s^2/3r}\right)n.$$
Since $f$ was arbitrary, this proves the theorem.
\end{proof}

The proof for general $k$ is, in this case, very similar. All that is
necessary is to throw out some `bad' vertices.

\begin{proof}[Proof of Theorem~\ref{jump}]
Let $n,r,s,k \in \N$, with $n \ge 16kr^2 + 4kr$, and let $f$ be an
$r$-colouring of the edges of $K_n$. If $s = 1$ then the result follows by
Mader's Theorem (since $n \ge 4kr + 1$), and the fact that $e^{-x} > 1 - x$
if $x > 0$, so assume that $s \ge 2$. Let $t \in [s - 1]$ (we shall again
eventually set $t = \lceil s/2 \rceil$, but we again delay making this
choice to emphasize the similarities with the previous proof), and let $G$
be a $k$-connected subgraph of $K_n$, using at most $t$ colours, of maximum
order. Let $V = V(K_n)$, $A = V(G)$, $B = V \setminus A$, and $T = f(E(G))$,
the set of colours used by $G$. Thus (assuming $|A| \le n - rk$), $|T| = t$.
By Mader's Theorem, we have $|A| \ge \ds\frac{n}{4r} \ge 4kr + k$.

Now, suppose there are at least $2kr\ds{{r-t} \choose t}$ vertices in $B$
which send at least $k$ edges of no more than $t$ colours into $A$. To be
more precise, given $v \in B$, let $$L_k(v) = \{\ell \in [r] : |\{u \in A :
f(uv) = \ell\}| \ge k\},$$ let $D = \{v \in B : |L_k(v)| \le t\}$, and
suppose that $|D| \ge 2kr\ds{{r-t} \choose t}$. Note that by
Observation~\ref{addvtx}, since $G$ is maximal, $L_k(v) \cap T = \emptyset$
for every $v \in B$. Thus, by the pigeonhole principle, there exists a
subset $S \subset [r] \setminus T$ of size $t$, and a subset $C \subset B$
of size $2kr$ such that $L_k(v) \subset S$ for every $v \in C$.

Consider the bipartite graph $H$, with parts $A$ and $C$, and edges with
colours from $S$. Note that, by the definition of $C$, each vertex of $C$
sends at most $kr$ edges with colours from $\ol{S}$ into $A$, so $d_H(v) \ge
|A| - kr$ for every $v \in C$. Let $a = 2kr$ and $b = kr$, and recall that
$|A| \ge 4kr + k = 4b + k$, and that $|C| \ge a \ge 2k$.

We apply Lemma~\ref{21ktech} to $H$, with $a = 2kr$ and $b = kr$, to obtain
a $k$-connected subgraph of $H$ on at least $$|A| + |C| - \frac{2k^2r^2}{2kr
- k + 1} \; > \; |A| + |C| - 2kr \; \ge \; |A|$$ vertices. This subgraph
uses at most $t$ colours, and so this contradicts the maximality of $G$.
Thus $|D| \le 2kr\ds{{r-t} \choose t}$.

Let $B' = B \setminus D$, so each vertex in $B'$ sends at least $k$ edges of
at least $t + 1$ different colours into $A$, i.e., $|L_k(v)| \ge t + 1$ for
every $v \in B'$. For each vertex $v \in B'$, choose a list $L(v) \subset
L_k(v)$ of size $t + 1$. So for each $v \in B'$ and $\ell \in L(v)$, there
exist at least $k$ vertices $u \in A$ such that $f(uv) = \ell$.

The remainder of the proof now goes through exactly as before, since by
Observation~\ref{addvtx}, for each vertex $v \in B'$ the vertices $A \cup
\{v\}$ are $k$-connected by the colours $T \cup \{\ell \}$ if $\ell \in
L(v)$. The reader who feels comfortable with this fact may therefore safely
`jump' to the end of the proof. For the remaining readers, and for
completeness, we shall repeat the argument.

So let $\eps > 0$, and let $\T = \{S \subset [r] \setminus T : |S| = s -
t\}$. Suppose that $m(n,r,s,k) \le n - \eps |B'|$. This means that for every
set $S \in \T$, the largest $k$-connected component in $K_n$, using only the
colours $S \cup T$, and containing $G$, avoids at least $\eps |B'|$ vertices
of $B'$. Hence, for each $S \in \T$ there are at least $\eps |B'|$ vertices
$v \in B'$ such that $S \cap L(v) = \emptyset$, by Observation~\ref{addvtx}.
For each $S \in \T$, let $M(S) = \{v \in B' : S \cap L(v) = \emptyset\}$.

Now, observe that for each vertex $v \in B'$, there are exactly
$\ds{{r-2t-1} \choose {s-t}}$ sets $S \in \T$ with $v \in M(S)$. So, summing
over $\T$, we obtain
$$\eps |B'| \ds{{r-t} \choose {s-t}} \; \le \; \sum_{v \in B'} \sum_{s \in
\T} I[v \in M(S)] \; = \; |B'| \ds{{r-2t} \choose {s-t}},$$ as before, and
therefore $$\eps \; \le \; \left( \ds\frac{r - 2t - 1}{r - t} \right)^{s-t}
\; < \; \operatorname{exp} \left( \ds\frac{-(t + 1)(s - t)}{r - t}
\right).$$
Now, setting $t = \lceil s/2 \rceil$ to (approximately) maximize
$\ds\frac{(t + 1)(s - t)}{r - t}$, and noting that $\ds\frac{\left( \lceil
s/2 \rceil + 1 \right) \lfloor s/2 \rfloor}{r - \lceil s/2 \rceil} >
\ds\frac{s^2}{4r}$, we obtain
\begin{align*}
M(f,n,r,s,k) & \ge \; |A| + (1 - \eps)|B'| \; \ge \; |A| + (1 - \eps) \left(
|B| - 2kr\ds{r \choose \lceil s/2 \rceil} \right)\\
& \ge \; \left( 1 - e^{-s^2/4r} \right) n - 2kr\ds{r \choose \lceil s/2
\rceil}.
\end{align*}
Since $f$ was arbitrary, this proves the theorem.
\end{proof}

\begin{rmk}
Using induction, as in the proof of Theorem~\ref{rootk=1}, one can slightly
improve this bound.\\
\end{rmk}

\section{Further Problems}\label{kprobs}

There is a great deal about the function $m(n,r,s,k)$ that we do not know.
In this section we shall discuss some of the most obvious and intriguing of
these open questions. We begin with the following corollary of
Theorem~\ref{jump} and Lemma~\ref{henry}. It demonstrates the rather
embarrassing state of our knowledge in the range $2 < s \ll \sqrt{r}$.

\begin{cor}
There exist constants $C, C' \in \RR$ such that
$$Cs^2 \; \le \; \ds\frac{r}{n} \: m(n,r,s,k) \; \le \; C' \min \left\{ 2^s
, s\sqrt{r} \right\}$$ for every $r,s,k \in \N$ with $s^2 < r$, and $n$
sufficiently large.
\end{cor}

\noindent In particular, we do not know whether the function $$g(s) =
\ds\liminf_{k \to \infty} \ds\liminf_{r \to \infty} \left( r \ds\liminf_{n
\to \infty} \left( \frac{1}{n} \: m(n,r,s,k) \right) \right)$$ grows like a
polynomial or an exponential function (or something in between!). We
conjecture that the upper bound is correct in the range $s \ll \log(r)$.

\begin{conj}\label{conj2^s}
Let $2 \le s,k \in \N$ be fixed. If $r > 4^s$, and $n$ is sufficiently
large, then $$m(n,r,s,k) \ge \ds\frac{2^sn}{r+1} - O(k).$$
\end{conj}

\noindent We suspect that Conjecture~\ref{conj2^s} is not easy, and pose the
following much weaker statements as open problems.

\begin{prob}
Prove any of the following.
\begin{enumerate}
\item[$(i)$] $g(s) > (1 + \eps)^s$ for some $\eps > 0$ and every $s \in \N$.
\item[$(ii)$] $g(s) < (1 + \eps)^s$ for every $\eps > 0$ and sufficiently
large $s$.
\item[$(iii)$] $g(s) = O(s^t)$ for some $t \in \N$.
\item[$(iv)$] $g(s) = \Omega(s^t)$ for every $t \in \N$.
\end{enumerate}
\end{prob}

\noindent When $r \gg s\sqrt{r} \gg 2^s$, we suspect that the upper bound in
Theorem~\ref{jump} becomes optimal, and $m(n,r,s,k) =
\Theta\left(\ds\ds\frac{sn}{\sqrt{r}}\right)$, but at present we seem a long
way from proving such a result.

We proved that the function $m(n,r,s,k)$ is `small' when $s \ll \sqrt{r}$
and `big' when $s \gg \sqrt{r}$. But what about when $s = \Theta(\sqrt{r})$?
What is the exact nature of this phase change? Theorem~\ref{jump} gives us
(roughly) the bounds $$\left(1 - e^{c/4}\right)n \; \le \; m(n,r,\lfloor
c\sqrt{r} \rfloor ,k) \; \le \; \ds\frac{cn}{\sqrt{2}}$$ when $n$ is
sufficiently large compared to $r$. Again we conjecture that the upper bound
is correct.

\begin{conj}\label{rootrconj}
Let $c \in (0,\sqrt{2}]$. Then $$h(c) \; = \; \liminf_{k \to \infty}
\liminf_{r \to \infty} \liminf_{n \to \infty} \left( \frac{1}{n} m(n, r,
\lfloor c\sqrt{r} \rfloor, k) \right) \; = \; \ds\frac{c}{\sqrt{2}}.$$
\end{conj}

Although we would really like to determine $h(c)$ exactly, we would in fact
be very happy with an answer to either of the following, more basic
questions.

\begin{qu}
Does there exist a constant $c \in \RR$ such that $h(c) = 1$?
\end{qu}

\begin{qu}
Is $\ds\lim_{c \to 0} \ds\frac{h(c)}{c} > 0$?
\end{qu}

Finally, we have a question about the phase transition at $2s = r$. We would
like to know the value of $m(n,2s-1,s,k)$; in other words, what does the
function jump \emph{to} when $r$ is odd? For $s = 2$ we showed that the
answer is $n - k + 1$, and it is tempting to guess that this is always the
correct answer, but we believe this to be false. More precisely we make the
following conjecture. The rather strange right-hand side is derived from a
fairly complicated construction, which we found and then lost! Since we
cannot prove it, we state it as a conjecture.

\begin{conj}
Let $n,s,k \in \N$, with $s \ge 3$ and $n$ and $k$ sufficiently large. Then
$$m(n,2s-1,s,k) \le n - \left( \frac{2s(r-2)}{r^2 - 2r + 2} \right)k = n -
\left( 1 + \frac{2s - 5}{4(s - 1)^2 + 1} \right)k,$$ where $r = 2s - 1$.
\end{conj}

\begin{prob}
Determine the value of $m(n,2s-1,s,k)$ for every $s,k \in \N$ and all
sufficiently large $n$.
\end{prob}

\section{Acknowledgements}

The authors would like to thank B\'ela Bollob\'as for suggesting the problem
to them, and for his ideas and encouragement. They would also like to thank
Microsoft Research, where part of this research was carried out.

\end{document}